\documentclass[leqno]{amsart}
\usepackage{amsfonts,amssymb,amsmath,amsgen,amsthm}
\usepackage{hyperref}

\theoremstyle{plain}
\newtheorem{theorem}{Theorem}[section]
\newtheorem{definition}[theorem]{Definition}
\newtheorem{assumption}[theorem]{Assumption}
\newtheorem{lemma}[theorem]{Lemma}
\newtheorem{corollary}[theorem]{Corollary}
\newtheorem{proposition}[theorem]{Proposition}
\theoremstyle{remark}
\newtheorem{remark}[theorem]{Remark}

\newtheorem{example}[theorem]{Example}

\def\C{{\mathbf C}}
\def\R{{\mathbf R}}
\def\N{{\mathbf N}}
\def\Sch{{\mathcal S}}
\def\F{\mathcal F}

\def\virgp{\raise 2pt\hbox{,}}

\def\({\left(}
\def\){\right)}
\def\<{\left\langle}
\def\>{\right\rangle}
\def\le{\leqslant}
\def\ge{\geqslant}

\def\Eq#1#2{\mathop{\sim}\limits_{#1\rightarrow#2}}
\def\Tend#1#2{\mathop{\longrightarrow}\limits_{#1\rightarrow#2}}

\def\d{{\partial}}
\def\l{\lambda}
\def\om{\omega}
\def\Om{\Omega}
\def\si{{\sigma}}
\def\eps{\varepsilon}
\DeclareMathOperator{\RE}{Re}
\DeclareMathOperator{\IM}{Im}

\numberwithin{equation}{section}

\begin{document}

\title[NLS with time dependent potential]{Nonlinear
  Schr\"odinger equation with time dependent potential}
\author[R. Carles]{R\'emi Carles}
\address{Univ. Montpellier~2\\Math\'ematiques
\\CC~051\\F-34095 Montpellier}
\address{CNRS, UMR 5149\\  F-34095 Montpellier\\ France}
\email{Remi.Carles@math.cnrs.fr}
\thanks{This work was supported by the French ANR project
  R.A.S. (ANR-08-JCJC-0124-01)}  

\begin{abstract}
  We prove a global well-posedness result for defocusing nonlinear
  Schr\"odinger equations with time dependent  potential. We then
  focus on  time dependent harmonic potentials. This
  aspect is motivated by 
  Physics (Bose--Einstein condensation), and appears also as a
  preparation for the analysis of the propagation of wave packets in
  a nonlinear context. The main aspect considered here is the growth
  of high Sobolev norms of the solution. 
\end{abstract}
\maketitle

\section{Introduction}
\label{sec:intro}
Let $d\ge 1$, and for $x\in \R^d$, consider the nonlinear
Schr\"odinger equation 
\begin{equation}
  \label{eq:NLSP}
  i\d_t u + \frac{1}{2}\Delta u = V(t,x) u + \lambda \lvert
  u\rvert^{2\si}u\quad ;\quad u_{\mid t=0}=u_0, 
\end{equation}
where $V\in\R$ is locally bounded in time and subquadratic in space,
$\l\in \R$, and the nonlinearity is energy-subcritical ($\si<2/(d-2)$
if $d\ge 3$). We prove that the solution exists and is global in  
\begin{equation*}
  \Sigma = \left\{ f \in H^1(\R^d)\quad ;\quad x\mapsto |x|f(x)\in
    L^2(\R^d)\right\}, 
\end{equation*}
provided that $u_0\in \Sigma$, $\si<2/d$ (mass-subcritical
nonlinearity), or $\si\ge 2/d$ and $\l\ge 0$ (defocusing 
nonlinearity). We then focus on the case where $V$ is exactly
quadratic in $x$:
\begin{equation}
  \label{eq:NLSharmo}
  i\d_t u + \frac{1}{2}\Delta u = \frac{1}{2}\sum_{j=1}^d\Om_j(t)x_j^2
  u + \lambda \lvert u\rvert^{2\si}u\quad ;\quad u_{\mid t=0}=u_0,
\end{equation}
where $\Om_j\in \R$, with $\Om_j\in C^1(\R)$. In the isotropic case
($\Om_j=\Om$ for all $j$), we show that the above result is optimal in
the sense that for all $\Om\in C(\R;\R)$, if
$\l<0$ and $\si=2/d$, there exist blow-up solutions. We also investigate the
growth of high Sobolev norms for large time. 
\smallbreak

There are at least two motivations to study \eqref{eq:NLSharmo}. In
Physics, this external potential may model a time dependent confining magnetic
potential:
\eqref{eq:NLSharmo} appears in Bose--Einstein condensation, typically
for $\si=1$ (or $\si=2$ sometimes in the one-dimensional case $d=1$),
see e.g. \cite{CD96,GRPGV01,RKRP08}. Equation~\eqref{eq:NLSharmo} also
appears as an envelope equation in the nonlinear propagation of wave
packets. In the linear case, consider
\begin{equation*}
  i\eps \d_t \psi^\eps +\frac{\eps^2}{2}\Delta \psi^\eps =
  V(x)\psi^\eps\quad ;\quad
  \psi^\eps(0,x)=\frac{1}{\eps^{d/4}}\varphi\(\frac{x-x_0}{\sqrt\eps}\)
  e^{i(x-x_0)\cdot \xi_0/\eps}. 
\end{equation*}
In the limit $\eps\to 0$, $\psi^\eps$ can be approximated as follows:
\begin{equation*}
  \psi^\eps(t,x)\Eq \eps 0 \frac{1}{\eps^{d/4}}
  u\(t,\frac{x-x(t)}{\sqrt\eps}\) e^{i\phi(t,x)/\eps},
\end{equation*}
where $(x(t),\xi(t))$ is given by the Hamiltonian flow associated to
$H= \frac{|\xi|^2}{2}+V(x)$, with initial data $(x_0,\xi_0)$, 
\begin{equation*}
  \phi(t,x)= \(x-x(t)\)\cdot \xi(t) +\int_0^t
  \(\frac{1}{2}|\xi(\tau)|^2-V\(x(\tau)\)\)d\tau, 
\end{equation*}
and $u$
is given by the equation
\begin{equation*}
  i\d_t u +\frac{1}{2}\Delta u = \frac{1}{2}\<Q(t)x,x\>u\quad ;\quad
  u_{\mid t=0}=\varphi.
\end{equation*}
Here, $Q$ is defined by $Q(t)={\rm Hess}V\(x(t)\)$, the Hessian of $V$
at point $x(t)$; see e.g. \cite{CR97}. We note that the external
potential in this case has the form presented in
\eqref{eq:NLSharmo}. To study the nonlinear propagation of wave 
packets, another parameter must be taken into account: the size of the
initial data. In  \cite{CF09ehrenfest}, it is shown that there exists a
critical size (depending on the nonlinearity and the space dimension),
corresponding to a certain power of $\eps$: for initial data which are
smaller (as $\eps\to 0$) than this critical size, then the
nonlinearity is negligible and we retrieve the same description as
above; for initial data which have the critical size, we have a
similar description, up to the fact that the envelope equation is now
nonlinear, of the form \eqref{eq:NLSharmo}. To study the propagation
of wave packets over large times (typically, up to -- an analogue of
-- Ehrenfest time), one has to understand the large time behavior of
solutions to \eqref{eq:NLSharmo}. 

\begin{remark}[Time dependent nonlinearity]
With little modification, we could also consider the more general equation
  \begin{equation}
    \label{eq:NLnonauto}
    i\d_t u + \frac{1}{2}\Delta u =
    \frac{1}{2}\sum_{j=1}^d\Om_j(t)x_j^2 
  u + h(t) \lvert u\rvert^{2\si}u\quad ;\quad u_{\mid t=0}=u_0,
  \end{equation}
where $h\in C^\infty(\R;\R)$. Following \cite{CW92} (see also
\cite{CazCourant}), the regularity assumption on $h$ could be
weakened. We choose to consider an autonomous nonlinearity in most of
this paper though.
\end{remark}

\begin{remark}[Complete integrability]
 The cubic one-dimensional case $d=\si=1$ is
special: if  $\Om=0$, 
then the equation is completely integrable
(\cite{AblowitzClarkson}). More generally,  \eqref{eq:NLnonauto}
has a Lax pair (recall that  $d=\si=1$) provided that  $\Om$ and $h$
are related through the identities (\cite{SHB07,kha09}):
\begin{equation*}
  \Om(t) = \ddot f(t)-\dot f(t)^2\quad ;
\quad h(t) = ae^{f(t)}\quad ;\quad a\in \R,\quad f\in C^\infty(\R;\R). 
\end{equation*}
The case where the above
relation is not satisfied is included in Proposition~\ref{prop:L2}. 
\end{remark}

The assumption we make on the external potential $V$ is inspired by
\cite{Fujiwara79}:
\begin{assumption}\label{hyp:V}
  $V\in L^\infty_{\rm loc}(\R_t\times \R_x^d)$ is smooth with respect
  to the space variable: for (almost) all $t\in \R$, $x\mapsto V(t,x)$
  is a $C^\infty$ map. Moreover, it is subquadratic in space: 
  \begin{equation*}
    \forall T>0, \ \forall \alpha \in \N^d, \ |\alpha|\ge 2, \quad
    \d_x^\alpha V\in L^\infty([-T,T]\times \R^d). 
  \end{equation*}
\end{assumption}
Note that this assumption does not involve spetral properties of
$V$, and are very little demanding concerning the dynamical properties
pour the associated Hamiltonian. The time dependent harmonic potential
that we consider in \eqref{eq:NLSharmo} is of course a peculiar case
of such potentials $V$.

\subsection{$L^2$-subcritical case}
\label{sec:L2sub}

When the energy is $L^2$-subcritical ($\si<2/d$), we have:
\begin{proposition}\label{prop:L2}
 Let $\l\in \R$, $0<\si<2/d$ and $V$ satisfying
 Assumption~\ref{hyp:V}. For all $u_0\in L^2(\R^d)$, 
  \eqref{eq:NLSP} has a unique solution
  \begin{equation*}
    u\in C\(\R;L^2(\R^d)\)\cap L^{\frac{4\si+4}{d\si}}_{\rm
      loc}\(\R;L^{2\si+2}(\R^d)\). 
  \end{equation*}
Moreover, its $L^2$-norm is independent of time:
\begin{equation*}
  \|u(t)\|_{L^2(\R^d)}=\|u_0\|_{L^2(\R^d)} ,\quad \forall t\in \R. 
\end{equation*}
\end{proposition}
\begin{proof}[Sketch of the proof]
  In view of \cite{Fujiwara79,Fujiwara}, local in time Strichartz
  estimates are available. Therefore, one can reproduce the original proof of
  \cite{TsutsumiL2} (see also \cite{CazCourant,TaoDisp}), in order to
  infer the result. 
\end{proof}

\subsection{Energy subcritical case}

In order to encompass the physical case $\si=1$ when $d=2$ or $3$, we
need to consider the case $\si\ge 2/d$. We shall restrict our
attention to $H^1$-subcritical nonlinearities: $\si<2/(d-2)$ when
$d\ge 3$. To solve \eqref{eq:NLSharmo}, even locally in time, one needs to work
in $\Sigma$, and not only $H^1$: symmetry is needed on physical and
frequency sides, unless $V$ is sublinear
\cite{CaPM08}. Local existence in $\Sigma$ then follows from the
dispersive estimates in \cite{Fujiwara79,Fujiwara}: one can work as
in the case $V\equiv 0$ (where it is possible to work in $H^1(\R^d)$
only). Instead of considering only $(u,\nabla u)$ as the unknown
function, one can consider $(u,\nabla u,xu)$: the three functions are related
through a closed family of estimates, and we get:
\begin{proposition}\label{prop:Sigmaloc}
  Let $\l\in \R$, $V$ satisfying
 Assumption~\ref{hyp:V},  and $\si>0$ with $\si<2/(d-2)$ if $d\ge 3$. For
  $u_0\in \Sigma$, there exists 
  $T$  and a
  unique solution $u$ solution to \eqref{eq:NLSP}, such that
  \begin{equation*}
    u,\nabla u,xu \in C\(]-T,T[;L^2(\R^d)\)\cap L^{\frac{4\si+4}{d\si}}_{\rm
      loc}\(]-T,T[;L^{2\si+2}(\R^d)\). 
  \end{equation*}
Moreover, its $L^2$-norm is independent of time:
\begin{equation*}
  \|u(t)\|_{L^2(\R^d)}=\|u_0\|_{L^2(\R^d)} ,\quad \forall t\in ]-T,T[. 
\end{equation*}
\end{proposition}
Since in \cite{Fujiwara79,Fujiwara}, only bounded time intervals are
considered, we give more precisions about this result in
\S\ref{sec:global} in order to consider global in time solutions: we
have been careful in the statement of Proposition~\ref{prop:Sigmaloc}
not to write that $T$ depends only on $\|u_0\|_\Sigma$. To infer
global existence results, we wish to replace the initial time $t=0$ in
\eqref{eq:NLSP} with $t=t_0\ge 0$: under the assumptions of
Proposition~\ref{prop:Sigmaloc}, it is not guaranteed that the
corresponding $T$ is independent of $t_0$. However, it is proved in
\cite{Fujiwara79,Fujiwara}  local dispersive estimates are available,
uniformly on finite time intervals.
\smallbreak

The natural candidate for an energy in the case of \eqref{eq:NLSP}
is 
\begin{equation}\label{eq:energy}
  E(t)=\frac{1}{2}\|\nabla u(t)\|_{L^2}^2
  +\frac{\lambda}{\si+1}\|u(t)\|_{L^{2\si+2}}^{2\si+2} +
  \int_{\R^d} V(t,x) \lvert u(t,x)\rvert^2dx. 
\end{equation}
\begin{proposition}
  Under the assumptions of Proposition~\ref{prop:Sigmaloc}, if in
  addition $V$ is $C^1$ with respect to $t$, and $\d_t V$ satisfies
  Assumption~\ref{hyp:V}, then $E\in
  C^1(]-T,T[;\R)$, and its evolution is given by
  \begin{equation}\label{eq:evolenergy}
    \frac{dE}{dt} = \int_{\R^d} \d_t V(t,x) \lvert u(t,x)\rvert^2dx.
  \end{equation}
\end{proposition}
The proof of the above result is straightforward, by following the
same lines as in the justification of similar evolution laws in,
e.g., \cite{CazCourant}. 
\begin{theorem}[Global existence in $\Sigma$]\label{theo:global}
  Let $\l\in \R$,  $\si>0$ with $\si<2/(d-2)$ if $d\ge 3$, and $V$
  satisfying Assumption~\ref{hyp:V}. For
  $u_0\in \Sigma$, we can take $T=+\infty$ in
  Proposition~\ref{prop:Sigmaloc} in the following cases: 
  \begin{itemize}
  \item $\si<2/d$.
\item $\si\ge 2/d$ and $\l\ge 0$, provided $V$ is $C^1$ in $t$ and
  $\d_tV$ satisfies Assumption~\ref{hyp:V}. 
  \end{itemize}
\end{theorem}
\begin{remark}\label{rem:dcds}
  This result extends the main one in \cite{CaDCDS}, where typically the
  (time independent) potential $-\om_1^2x_1^2+\om^2_2 x_2^2$ is
  considered. It is 
  established in \cite{CaDCDS} that if $\l>0$ and $\om_1\gg 1+\om_2$,
  then the solution to \eqref{eq:NLSharmo} is global, and there is
  scattering. The present theorem extends the existence part, but of
  course the assumptions of Theorem~\ref{theo:global} are too general
  to expect a scattering result: in the case of \eqref{eq:NLSharmo}
  with $\Om_j=1\ \forall j$,
  for instance, one can
  construct periodic solutions to 
  \eqref{eq:NLSharmo}, of the form $u(t,x)= e^{-i\om
    t}\psi(x)$. Indeed, this amounts to finding a non-trivial solution
  to the elliptic problem
  \begin{equation}\label{eq:elliptic}
    \om \psi = H  \psi  +\l
    |\psi|^{2\si}\psi,\text{ where }H=-\frac{1}{2}\Delta+\frac{|x|^2}{2}.
  \end{equation}
Introduce the quantities
  \begin{align*}
    I(\psi)&= \frac{1}{2}\<H\psi,\psi\> -\frac{\om}{2}\<\psi,\psi\>,\\
 M&=\left\{ \psi\in \Sigma\ ; \ 
    \frac{1}{\si+1} \int_{\R^d}|\psi(x)|^{2\si+2}dx=1\right\}, 
  \end{align*}
and consider
\begin{equation*}
  \delta =\inf_{\psi\in M}I(\psi).
\end{equation*}
If $\om>d/2$ (the lowest eigenvalue of the harmonic oscillator), then
\eqref{eq:elliptic} has a non-trivial solution for $\l > 0$. If $\om<d/2$, then
\eqref{eq:elliptic} has a non-trivial solution for $\l < 0$. See
e.g. \cite{CaDPDE} for more details. 
\end{remark}
Theorem~\ref{theo:global} shows that in the usual cases where global
existence is known without a potential, the introduction of a smooth
subquadratic potential $V$ does not change this property, regardless of
the time dependence of $V$ with respect to time. 
In the case $\si=2/d$ and $\l<0$, we prove that finite time blow-up
does occur for time dependent potentials, like in the case with no
potential: 
\begin{proposition}[Finite time blow-up]\label{prop:blowup}
  Let $\si=2/d$ and $\l<0$. Consider \eqref{eq:NLSharmo} with an
  isotropic potential: $\Om_j=\Om\in C(\R;\R)$ is independent of
  $j$. There exist blow-up solutions for \eqref{eq:NLSharmo}: we can
  find $T>0$, and $u\in C(]0,T];\Sigma)$ solving \eqref{eq:NLSharmo}
  such that
  \begin{equation*}
    \|\nabla u(t)\|_{L^2}\Tend t 0 \infty. 
  \end{equation*}
\end{proposition}

\subsection{Growth of higher order Sobolev norms}
\label{sec:introgrowth}

We now stick to the case of time dependent harmonic potentials,
\eqref{eq:NLSharmo}. 
In view of the analysis of nonlinear wave packets in a semi-classical
regime (\cite{CF09ehrenfest}), the evolution of weighted Sobolev norms
of $u$ over 
large time intervals is needed. 
\smallbreak

Consider first the autonomous isotropic case: $\Om_j=\Om$ is a constant.

If $\Om=0$, then at least when $\l\ge 0$ and $2/d\le \si<2/(d-2)$,
the Sobolev norms of $u$ are bounded, $u\in L^\infty(\R;H^k(\R^d))$
(provided the nonlinearity is sufficiently smooth), since we know that
there is scattering in $H^1$ (because we know that there is scattering
in $\Sigma$, since scattering in $H^1$ only is not known so far in the
case $\si=2/d$); see
e.g. \cite{Wa04}. The momenta of $u$  
grow algebraically in time (see \cite{Be01}). We give a short
alternative proof of these properties in an appendix.
\smallbreak

If $\Om>0$, then $u\in L^\infty(\R;\Sigma)$, as proved
by \eqref{eq:evolenergy}. The existence of periodic solutions to the
nonlinear problem (see Remark~\ref{rem:dcds}) shows that we may also
have $u\in L^\infty(\R;H^k)$ and $|x|^ku \in L^\infty(\R;L^2)$ for all
$k\in \N$. 
\smallbreak

If $\Om<0$ (repulsive harmonic
  potential), then it is proved in \cite{CaSIMA} that every defocusing
  $H^1$-subcritical nonlinearity is short range as far as scattering
  theory is concerned: if $\l>0$ in \eqref{eq:NLSharmo}, then
  $u(t)\sim U(t)u_+$ as $t\to +\infty$, for some $u_+\in \Sigma$,
  where $U(t) =\exp\(-it(-\frac{1}{2}\Delta
  +\frac{\Om}{2}|x|^2)\)$. Assume $\Om=-1$. Using Mehler's formula,
  and a decomposition 
  of $U(t)$ of the form $U=MD\F M$ as in \cite{TsutsumiSigma}
  originally (for the case $\Om=0$), we notice
  \begin{equation*}
    U(t)u_+(x)\Eq t{+\infty} \frac{1}{\sinh t}
    \widehat u_+\(\frac{x}{\sinh t}\) e^{i\frac{\cosh
     t}{\sinh t}|x|^2}. 
  \end{equation*}
This shows that the $L^2$ norms of $\nabla U(t)u_+$ and $x U(t)u_+$
grow exponentially in time. By the results in \cite{CaSIMA}, so do the
$L^2$ norms of $\nabla u$ and $x u$. Note that at least in the linear
case $\l=0$, we see that the $H^k$-norms of $u$ grow like $e^{kt}$ as
$t$ goes to infinity. 

\begin{definition}
 Let $u\in C(\R;\Sigma)$ be a solution to \eqref{eq:NLSharmo}, and
 $k\in \N$. \\
$\bullet$ $(Alg)_k$ is satisfied if
  there exists $A$ such that for all  admissible pair $(p,q)$, 
  \begin{equation*}
    \forall \alpha, \beta\in \N^d,\ |\alpha|+|\beta|\le k, \quad
    \left\lVert x^\alpha \d_x^\beta u\right\rVert_{L^p([0,t];L^q)}\lesssim
     t^A. 
  \end{equation*}
$\bullet$ $(Exp)_k$ is satisfied if
  there exists $C$ such that for all  admissible pair $(p,q)$, 
  \begin{equation*}
    \forall \alpha, \beta\in \N^d,\ |\alpha|+|\beta|\le k, \quad
    \left\lVert x^\alpha \d_x^\beta u\right\rVert_{L^p([0,t];L^q)}\lesssim
     e^{Ct}. 
  \end{equation*}
\end{definition}
We wish to consider smooth energy-subcritical nonlinearities. Since we
study homogeneous nonlinearities, we have to assume: $d\le 3$,
$\sigma\in \N$, with $\si=1$ if $d= 3$. 
\begin{corollary}\label{cor:alg}
 Let $d\le 3$, $\l\ge 0$ and $\sigma\in \N$, with $\si=1$ if $d=
 3$. Let $k\in \N$, 
  \begin{equation*}
    u_0 \in H^k\(\R^d\),\text{ with }|x|^ku_0\in L^2(\R^d). 
  \end{equation*}
If $\Om_j \in
 C^1(\R;\R)$ is compactly supported for all $j$, then $u$ has the property
 $(Alg)_k$. 
\end{corollary}
\begin{proof}
  We may assume that $\operatorname{supp}\Om\subset [-M,M]$. From
  Theorem~\ref{theo:global}, $u_{\mid t=M}\in \Sigma$. It is easy to
  check that the higher regularity is conserved as well: this is
  rather straightforward, since we consider an energy-subcritical
  nonlinearity. The corollary
  then follows from the case $\Om=0$, where $(Alg)_k$ is satisfied, as
  recalled in the appendix.
\end{proof}
In the case where the dependence of $\Om_j$ with respect to time is
not specified, the evolution of the energy \eqref{eq:energy} yields no
exploitable information. Even in the $L^2$-subcritical case, we will
see that proving exponential control requires some work.
\begin{proposition}[Exponential growth]\label{prop:growth}
  Let $d\le 3$, $\sigma\in \N$, with $\si=1$ if $d= 3$, $\Om_j\in
  C(\R;\R)$ be locally Lipschitzean, $k\in \N$, $k\ge
  1$, and  
  \begin{equation*}
    u_0 \in H^k\(\R^d\),\text{ with }|x|^ku_0\in L^2(\R^d). 
  \end{equation*}
$(Exp)_k$ is satisfied (at least) in the following cases:
\begin{itemize}
\item $\si=d=1$ ($L^2$-subcritical nonlinearity), and $\Om$ is bounded. 
\item $\si\ge 2/d$, $\l\ge 0$, and $\Om_j=\Om\le 0$ is independent of $j$
  (isotropic repulsive  potential). 
\end{itemize}
\end{proposition}
\begin{remark}[Optimality]
   When the potential is repulsive and time-independent ($\Om=-1$
   typically), the 
   exponential growth is sharp, and $C$ does depend on $k$ ($C=k$
   when $\Om=-1$), as discussed above.  
\end{remark}

\begin{remark}
 We prove that in the case of an isotropic repulsive  potential, there
 is scattering provided $\si\ge 2/d$
 (Proposition~\ref{prop:scattrepuls}): morally, $(Exp)_k$ is
 satisfied because it is satisfied in the linear setting (case
 $\l=0$). However, this property on the linear solution demands a
 justification; the key is Lemma~\ref{lem:edorep}. 
\end{remark}

\subsection{Outline of the paper}
\label{sec:outline}

In Section~\ref{sec:global}, we prove Theorem~\ref{theo:global}. We
then focus on the study of Equation~\eqref{eq:NLSharmo}. In
Section~\ref{sec:mehler}, we derive a generalized Mehler formula to
express the fundamental solution associated to the linear equation,
\eqref{eq:NLSharmo} with $\l=0$. In Section~\ref{sec:lens}, we
generalize a lens transform, known in the case of isotropic
time-independent quadratic potentials, to the case of isotropic
time-dependent quadratic potentials. This allows us to infer
Proposition~\ref{prop:blowup}. In Section~\ref{sec:champs}, we
introduce some vector fields, which correspond to Heisenberg
derivatives, and yield interesting evolution laws when the potential
is isotropic. In Section~\ref{sec:growth}, we examine the large time
behavior of solutions to \eqref{eq:NLSharmo}, and prove
Proposition~\ref{prop:growth}. Finally, we show in an appendix that when
$V=0$, for large time, the solutions to \eqref{eq:NLSP} have bounded
Sobolev norms and algebraically growing momenta, provided there
is scattering.

\section{Global existence in $\Sigma$: proof of Theorem~\ref{theo:global}}
\label{sec:global}

\subsection{Strichartz estimates}
\label{sec:strichartz}

We first recall some results established in
\cite{Fujiwara79,Fujiwara}. Consider $V$ satisfying
Assumption~\ref{hyp:V}. It is  established in 
\cite{Fujiwara79} that one can define $U(t,s)$ as $u(t,x) =
U(t,s)\varphi(x)$, where 
\begin{equation}\label{eq:lins}
  i \d_t u+\frac{1}{2}\Delta u
  = V(t,x) u 
  \quad ;\quad u(s,x)=\varphi(x),
\end{equation}
along
with the following properties: 
\begin{itemize}
\item $U(t,t)=\rm{Id}$.
\item The map $(t,s)\mapsto U(t,s)$ is strongly continuous.
\item $U(t,s)^*=U(t,s)^{-1}$.
\item $U(t,\tau)U(\tau,s)=U(t,s)$. 
\item $U(t,s)$ is unitary on $L^2$:
  $\|U(t,s)\varphi\|_{L^2}= \|\varphi\|_{L^2}$. 
\end{itemize}
In addition, we know from
\cite{Fujiwara} that for all $T>0$, $t,s\in [-T,T]$, 
\begin{equation}\label{eq:disp}
\lVert U(t,0)U(s,0)^*\varphi\rVert_{L^\infty(\R^d)} = \lVert
  U(t,s)\varphi\rVert_{L^\infty(\R^d)}\le 
  \frac{C}{|t-s|^{d/2}} \|\varphi\|_{L^1(\R^d)}, 
\end{equation}
provided that $|t-s|<\eta$. It is implicitly assumed in
\cite{Fujiwara79} that $\eta$ may depend on $T$; in
Example~\ref{ex:StrichartzUSTC} below, we show that 
it does indeed, in the case of the time dependent harmonic potential,
if the functions $\Om_j$ are not bounded.   

Recall the standard definition in the context of Schr\"odinger
equations:
\begin{definition}\label{def:adm}
 A pair $(p,q)$ is admissible if $2\le q
  <\frac{2d}{d-2}$ ($2\le q\le\infty$ if $d=1$, $2\le q<
  \infty$ if $d=2$)  
  and 
$$\frac{2}{p}=\delta(q):= d\left( \frac{1}{2}-\frac{1}{q}\right).$$
\end{definition}
The general results on Strichartz estimates (see e.g. \cite{KT}) then
yield, as a consequence of the dispersive estimate \eqref{eq:disp}: 
\begin{proposition}\label{prop:strichartz}
  Recall that $U(t,s)$ is defined by \eqref{eq:lins}, where
$V$ satisfies Assumption~\ref{hyp:V}. Let $T>0$.  There exists
$\eta>0$ such that the following holds:\\ 
$(1)$ For any admissible pair $(p,q)$, there exists $C_{q}$  such that
$$
\|U(\cdot,s) \varphi\|_{L^{p}([s,s+\eta];L^{q})} \le C_q 
\|\varphi \|_{L^2},\quad \forall \varphi\in L^2(\R^d),\quad \forall
s\in ]-T,T-\eta[.
$$
$(2)$ 
For $s\in \R$, denote
\begin{equation*}
  D_s(F)(t,x) = \int_{s}^t U(t,\tau)F(\tau,x)d\tau. 
\end{equation*}
For all admissible pairs $(p_1,q_1)$ and~$
    (p_2,q_2)$,  there exists $C=C_{q_1,q_2}$ independent of $s\in ]-T,T-\eta[$
    such that  
\begin{equation}\label{eq:strichnl}
      \left\lVert D_s(F)
      \right\rVert_{L^{p_1}([s,s+\delta];L^{q_1})}\le C \left\lVert
      F\right\rVert_{L^{p'_2}\([s,s+\delta];L^{q'_2}\)},
\end{equation}
for all $F\in L^{p'_2}(I;L^{q'_2})$ and $0\le \delta\le \eta$.
\end{proposition}
\begin{example}[Standard harmonic oscillator]
  Assume that $V(t,x)=\frac{|x|^2}{2}$. The above result is then
  standard (see e.g. \cite{CazCourant}). The fact that one has to consider
  finite time intervals for the above result to be valid stems for
  instance from the existence of eigenvalues for the harmonic
  oscillator: let $g(x)= e^{-|x|^2/2}$ be the ground state associated
  to the harmonic potential, and denote $u(t,x)=e^{-itd/2}g(x)$. It
  solves
  \begin{equation*}
    i\d_t u +\frac{1}{2}\Delta u =\frac{|x|^2}{2}u\quad ;\quad u_{\mid
      t=0}=g. 
  \end{equation*}
We compute $\|u\|_{L^p(I;L^q)}= |I|^{1/p}\|g\|_{L^q}$, which shows
that Proposition~\ref{prop:strichartz} becomes false with
$\eta=\infty$. 
\end{example}
\begin{example}\label{ex:StrichartzUSTC}
  We show that in general, the above result is false with $T=\infty$. 
Let 
\begin{equation*}
  V(t,x)=\frac{1}{2}\Om(t)|x|^2. 
\end{equation*}
If  $\Om$ is not bounded, then the above
  uniformity with respect to $s$ fails: let
  \begin{equation*}
    \Om(t) = n^2\text{ if }4n+1=t_n \le t\le 4n+2. 
  \end{equation*}
Since we have
\begin{equation*}
  \(-\frac{1}{2}\Delta +\frac{n^2}{2}|x|^2\) e^{-n|x|^2/2}=
  \frac{nd}{2} e^{-n|x|^2/2},
\end{equation*}
the function $u(t,x)=e^{-ind(t-t_n)/2 -n|x|^2/2}$ solves \eqref{eq:lins} with
$s=t_n$. If Proposition~\ref{prop:strichartz} was true with
$T=\infty$, we would have: 
\begin{equation*}
  \|u\|_{L^p([4n+1,4n+1+\eta];L^q)}=
  \eta^{1/p}\(\frac{2\pi}{nq}\)^{d/(2q)}\le C \|u(t_n)\|_{L^2} = C
  \(\frac{\pi}{n}\)^{d/4} ,
\end{equation*}
where $C$ does not depend on $n$. For all $q>2$, letting $n$ go to
infinity leads to a contradiction. Since \eqref{eq:disp} implies
Proposition~\ref{prop:strichartz}, this shows \eqref{eq:disp} is valid
for $|t-s|<\eta$ where $\eta$ depends on $T$, unless $\Om$ is
bounded. 
\end{example}

\subsection{Local existence in $\Sigma$}
\label{sec:locSigma}

Since \eqref{eq:NLSP} is not autonomous, we consider the same
problem with a varying initial time:
\begin{equation}
  \label{eq:NLSP2}
  i\d_t u + \frac{1}{2}\Delta u = V(t,x)
  u + \lambda \lvert u\rvert^{2\si}u\quad ;\quad u_{\mid t=s}=u_0,
\end{equation}
with $s\in \R$. 
\begin{proposition}\label{prop:Sigmaloc2}
  Let $\l\in \R$,    $\si>0$
  with $\si<2/(d-2)$ if 
  $d\ge 3$, and $V$ satisfying Assumption~\ref{hyp:V}. Let $M>0$, and
  $s\in ]-M,M[$.  For  all $u_0\in \Sigma$, 
  there exists  $T=T(\|u_0\|_\Sigma,M)$  and a
  unique solution $u$ solution to \eqref{eq:NLSP2}, such that
  \begin{equation*}
    u,\nabla u,xu \in C\(]s-T,s+T[;L^2(\R^d)\)\cap
    L^{\frac{4\si+4}{d\si}}_{\rm 
      loc}\(]s-T,s+T[;L^{2\si+2}(\R^d)\). 
  \end{equation*}
Moreover, its $L^2$-norm is independent of time:
\begin{equation*}
  \|u(t)\|_{L^2(\R^d)}=\|u_0\|_{L^2(\R^d)} ,\quad \forall t\in ]s-T,s+T[. 
\end{equation*}
If $V$ is $C^1$ in $t$, then the  energy $E$ (defined by
\eqref{eq:energy}) evolves according to \eqref{eq:evolenergy}. 
\end{proposition}
\begin{proof}[Sketch of the proof]
  We present here only the main steps of the classical
  argument. Duhamel's formulation for \eqref{eq:NLSP2} reads
  \begin{equation*}
    u(t)=U(t,s)u_0 -i\l \int_s^t U(t,\tau)\(|u|^{2\si}u\)(\tau)d\tau. 
  \end{equation*}
Denote the right hand side by $\Phi^s(u)(t)$. 
Proposition~\ref{prop:Sigmaloc2} follows from a fixed point argument
in the space
\begin{equation*}
  X_T = \left\{u\in C(I_T;\Sigma)\ ;\ u,xu,\nabla u \in
    L^{\frac{4\si+4}{d\si}}\(I_T;L^{2\si+2}(\R^d)\) \right\},
\end{equation*}
where $I_T=]s-T,s+T[$.  
Introduce the Lebesgue exponents 
\begin{equation*}
  q=2\si+2\quad ;\quad p=\frac{4\si+4}{d\si}\quad ;\quad
  \theta=\frac{2\si(2\si+2)}{2-(d-2)\si}. 
\end{equation*}
Then $(p,q)$ is admissible, and
\begin{equation*}
  \frac{1}{q'}=\frac{2\si}{q}+\frac{1}{q}\quad ;\quad
  \frac{1}{p'}=\frac{2\si}{\theta} +\frac{1}{p}. 
\end{equation*}
Proposition~\ref{prop:strichartz} and H\"older inequality yield
\begin{align*}
  \|\Phi^s(u)\|_{L^p(I_T;L^q)\cap L^\infty (I_T;L^2)}&\le C \|u_0\|_{L^2}
  +C \left\lVert
    \lvert u\rvert^{2\si}u\right\rVert_{L^{p'}(I_T;L^{q'})}\\
&\le C \|u_0\|_{L^2}
  +C \|u\|_{L^\theta(I_T;L^q)}^{2\si} \|u\|_{L^p(I_T;L^q)}, 
\end{align*}
where $C$ is independent of $s\in [-M,M]$ and $T\le \eta$. Using Sobolev
embedding, 
\begin{equation*}
  \|\Phi^s(u)\|_{L^p(I_T;L^q)\cap L^\infty (I_T;L^2)} \le C \|u_0\|_{L^2}
  +C T^{2\si/\theta}\|u\|_{L^\infty(I_T;H^1)}^{2\si} \|u\|_{L^p(I_T;L^q)}.
\end{equation*}
We have
\begin{align*}
  \nabla \Phi^s(u)(t) &= U(t,s)\nabla u_0 -i\l \int_s^t
  U(t,\tau)\nabla \(|u|^{2\si}u\)(\tau)d\tau \\
&\quad -i\int_s^t U(t,\tau)\(\Phi^s(u)(\tau)\nabla V(\tau)\)d\tau. 
\end{align*}
We estimate the second term of the right hand side as above, and get,
since $\nabla V$ is sublinear by assumption:
\begin{align*}
  \|\nabla\Phi^s(u)\|_{L^p(I_T;L^q)\cap L^\infty (I_T;L^2)} &\le C
  \|\nabla u_0\|_{L^2} + C T^{2\si/\theta}\|u\|_{L^\infty(I_T;H^1)}^{2\si} \|\nabla
  u\|_{L^p(I_T;L^q)} \\
&\quad + C\|\Phi^s(u)\nabla V\|_{L^1(I_T;L^2)}\\
&\le C
  \|\nabla u_0\|_{L^2} + C T^{2\si/\theta}\|u\|_{L^\infty(I_T;H^1)}^{2\si} \|\nabla
  u\|_{L^p(I_T;L^q)} \\
&\quad + C T\|x\Phi^s(u)\|_{L^\infty(I_T;L^2)}+ C
T\|\Phi^s(u)\|_{L^\infty(I_T;L^2)}, 
\end{align*}
where, again, $C$ does not depend on $s\in [-M,M]$. We have similarly
\begin{equation*}
  \begin{aligned}
       \|x\Phi^s(u)\|_{L^p(I_T;L^q)\cap L^\infty (I_T;L^2)} 
&\le C
  \|x u_0\|_{L^2} + C T^{2\si/\theta}\|u\|_{L^\infty(I_T;H^1)}^{2\si} \|x
  u\|_{L^p(I_T;L^q)} \\
&\quad + C T\|\nabla \Phi^s(u)\|_{L^\infty(I_T;L^2)}.
  \end{aligned}
\end{equation*}
Choosing $T$ sufficiently small, one can then prove that $\Phi^s$ maps
a suitable ball in $X_T$ into itself. Contraction for the norm
$\|\cdot\|_{L^p(I_T;L^q)}$ is proved similarly, and one concludes by
remarking that $X_T$ equipped with this norm 
is complete. 
\end{proof}
We can now infer the analogue to the standard result (which is not
straightforward since we consider a non-autonomous equation, in the
presence of an external potential):
\begin{corollary}\label{cor:alternative}
  Let $\l\in \R$, $\si>0$
  with $\si<2/(d-2)$ if  $d\ge 3$,  $V$ satisfying
  Assumption~\ref{hyp:V}, and $u_0\in \Sigma$. Either the solution to 
  \eqref{eq:NLSharmo} is global in time (in the future), 
  \begin{equation*}
    u,\nabla u,xu \in C\(\R_+;L^2(\R^d)\)\cap
    L^{\frac{4\si+4}{d\si}}_{\rm 
      loc}\(\R_+;L^{2\si+2}(\R^d)\), 
  \end{equation*}
or there exists $T>0$, such that 
\begin{equation*}
  \|\nabla u(t)\|_{L^2}
\mathop{\longrightarrow}\limits_{t{\mathop{\rightarrow}\limits_<}  T} 
+\infty.  
\end{equation*}
\end{corollary}
\begin{proof}
Let $M>0$. 
  Proposition~\ref{prop:Sigmaloc2} shows that the only obstruction to
  well-posedness on  $[0,M]$ is the existence of a
  time $0<T<M$ such that 
\begin{equation*}
  \|x u(t)\|_{L^2}+\|\nabla u(t)\|_{L^2}
\mathop{\longrightarrow}\limits_{t{\mathop{\rightarrow}\limits_<}  T} 
 +\infty.  
\end{equation*}
So long as $u\in
C([0,t];\Sigma)$, we have (see e.g. \cite{CazCourant} for the
arguments that make the computation rigorous) 
\begin{equation}
  \label{eq:evolmoment}
   \frac{d}{dt}\int_{\R^d}x_j^2 \lvert u(t,x)\rvert^2dx = 2\IM
    \int_{\R^d}x_j \overline 
    u(t,x)\d_j u (t,x)dx.  
\end{equation}
Suppose  $u\in L^\infty ([0,T];H^1)$. Then the above formula,
Cauchy--Schwarz inequality and Gronwall lemma show that $x u\in L^\infty
([0,T];L^2)$,
hence a contradiction. The corollary follows, since $M>0$ is arbitrary. 
\end{proof}
Therefore, to prove global existence in $\Sigma$ in the
$H^1$-subcritical case, it suffices to exhibit \emph{a priori} bounds
for $\nabla u$  in $L^2$.

\subsection{$L^2$-subcritical case}
\label{sec:L2sub2}

In the case $\si<2/d$, recall that the classical argument of \cite{TsutsumiL2}
can be applied directly, to infer Proposition~\ref{prop:L2}. The
\emph{a priori} bound for $(\nabla u,xu)$ in $L^2$ then follows, by
resuming the computations presented in the proof of
Proposition~\ref{prop:Sigmaloc2}. Keeping the same notations, we have
in particular 
\begin{align*}
 & \|\nabla\Phi^s(u)\|_{L^p(I_T;L^q)\cap L^\infty (I_T;L^2)}
  +\|x\Phi^s(u)\|_{L^p(I_T;L^q)\cap L^\infty (I_T;L^2)}  \\
&\le C 
  \| u_0\|_{\Sigma} + C \|u\|_{L^\theta(I_T;L^q)}^{2\si} \(\|\nabla
  u\|_{L^p(I_T;L^q)}+\|x
  u\|_{L^p(I_T;L^q)}\) \\
&\quad + C T\(\|\Phi^s(u)\|_{L^\infty(I_T;L^2)}+\|x\Phi^s(u)\|_{L^\infty(I_T;L^2)}
+\|\nabla\Phi^s(u)\|_{L^\infty(I_T;L^2)}\) , 
\end{align*}
where we recall that
\begin{equation*}
    q=2\si+2\quad ;\quad p=\frac{4\si+4}{d\si}\quad ;\quad
  \theta=\frac{2\si(2\si+2)}{2-(d-2)\si},
\end{equation*}
and, in view of Proposition~\ref{prop:Sigmaloc2}, we know that $u=
\Phi^s(u)$. In the case $\si<2/d$, we have $1/p<1/\theta$, and thus
\begin{equation*}
  \|u\|_{L^\theta(I_T;L^q)} \le (2T)^{1/\theta-1/p}\|u\|_{L^p(I_T;L^q)}=
(2T)^{\frac{(2-d\si)(\si+1)}{2\si(2\si+2)}} \|u\|_{L^p(I_T;L^q)}.
\end{equation*}
By Proposition~\ref{prop:L2}, $u\in L^p_{\rm
  loc}(\R;L^q(\R^d))$. Splitting any given time interval $[-M,M]$ into
finitely many (tiny) pieces, we obtain an \emph{a priori} bound for $(\nabla
u,xu)$ in $L^\infty([-M,M];L^2)$. Since $M>0$ is arbitrary,
Corollary~\ref{cor:alternative}  
yields the first point of Theorem~\ref{theo:global}. 
\subsection{Defocusing energy-subcritical case}

We now consider the case $\l\ge 0$, with $\si<2/(d-2)$ if $d\ge 3$. 
To complete the proof of Theorem~\ref{theo:global}, we resume the
computation initiated in the proof of Corollary~\ref{cor:alternative},
in order to infer a virial identity:
\begin{lemma}\label{lem:viriel}
  Let $\l\in \R$,    $\si>0$ with $\si<2/(d-2)$ if
  $d\ge 3$, and $V$ satisfying Assumption~\ref{hyp:V}. Let $u_0\in
  \Sigma$, and $u\in C(]-T,T[;\Sigma)$ be the 
  solution to \eqref{eq:NLSP} given by
  Proposition~\ref{prop:Sigmaloc2} (case $s=0$). Denote
  \begin{equation*}
    y(t)= \int_{\R^d}\lvert x\rvert^2 \lvert u(t,x)\rvert^2dx. 
  \end{equation*}
Then $y\in C^2(]-T,T[)$, and satisfies
\begin{equation*}
  \frac{d^2y}{dt^2} = 2\|\nabla u(t)\|_{L^2}^2 -2\int_{\R^d}x\cdot
  \nabla V(t,x)
  \lvert u(t,x)\rvert^2dx+ 2\l \frac{ 
    d\si}{\si+1}\|u(t)\|_{L^{2\si+2}}^{2\si+2}   .
\end{equation*}
\end{lemma}
\begin{proof}
  We present the formal part of the proof, and refer to
  \cite{CazCourant} for the arguments that make the proof rigorous. We
  first resume the computation made in the course of the proof of
  Corollary~\ref{cor:alternative}.
Differentiating \eqref{eq:evolmoment} with respect to time again, we have:
\begin{align*}
 \frac{d^2}{dt^2}\|x_j u\|_{L^2}^2 &= 
2 \IM
    \int_{\R^d}x_j \d_t \overline 
    u\d_j u   + 2 \IM
    \int_{\R^d}x_j \overline 
     u\d_j \d_t u \\
&= -2  \IM
    \int_{\R^d}\(\overline 
    u +2 x_j \d_j \overline 
    u\)  \d_t u = 2\RE \int_{\R^d}\(\overline 
    u +2 x_j \d_j \overline 
    u\)  i\d_t u \\
&= 2\RE \int_{\R^d}\(\overline 
    u +2 x_j \d_j \overline 
    u\)  \( -\frac{1}{2}\Delta u +V(t,x) u +\l
    |u|^{2\si}u\) 
\end{align*}
The terms in factor of  $\overline   u$ simplify easily, and we infer:
\begin{align*}
 \frac{d^2}{dt^2}\|x_j u\|_{L^2}^2 
&= \|\nabla u\|_{L^2}^2 +2\int_{\R^d}V(t,x)|u(t,x)|^2dx + 2\l
\|u\|_{L^{2\si+2}}^{2\si+2} \\
&\quad -2\RE \int_{\R^d} x_j \d_j\overline u \Delta u +4 
 \RE \int_{\R^d} V(t,x)x_j  u\d_j \overline u \\
&\quad +4\l \RE
\int_{\R^d} x_j |u|^{2\si}u\d_j \overline u\\
&= \|\nabla u\|_{L^2}^2 + 2\int_{\R^d}V(t,x)|u(t,x)|^2dx+ 2\l
\|u\|_{L^{2\si+2}}^{2\si+2} \\
&\quad -\|\nabla u\|_{L^2}^2+ 2\|\d_j u\|_{L^2}^2 +2 \int_{\R^d}x_j
V(t,x)\d_j\(|u|^2\) \\
&\quad -\frac{2\l}{\si+1} \|u\|_{L^{2\si+2}}^{2\si+2}\\
&= 2\|\d_j u\|_{L^2}^2-2\int_{\R^d}x_j\d_j V(t,x)|u(t,x)|^2dx+2\l
\frac{\si}{\si+1} 
\|u\|_{L^{2\si+2}}^{2\si+2}. 
\end{align*}
The result then follows by summing over $j$. 
\end{proof}
To complete the proof of Theorem~\ref{theo:global}, fix $M>0$, and for
$t\in [0,M]$, let 
\begin{equation*}
  f(t) = y(t) + \lvert \dot y(t)\rvert. 
\end{equation*}
We have
\begin{equation*}
  \dot f(t)  \le  \lvert \dot y(t)\rvert + \lvert \ddot y(t)\rvert
 \le  \lvert \dot y(t)\rvert + 2 \|\nabla u\|_{L^2}^2 + C+ C
 y(t) + C \|u\|_{L^{2\si+2}}^{2\si+2},
\end{equation*}
where we have used Lemma~\ref{lem:viriel}, the estimate
\begin{equation*}
  \lvert x\cdot \nabla V(t,x)\rvert \le C(M)\(1+|x|^2\),\quad \forall
  (t,x)\in [0,M]\times\R^d, 
\end{equation*}
and the conservation of mass. 
Since $u_0\in\Sigma$, \eqref{eq:energy}--\eqref{eq:evolenergy}
(this is where we have to assume that $V$ is $C^1$ in $t$) yield
\begin{equation*}
  \|\nabla u\|_{L^2}^2+  \|u\|_{L^{2\si+2}}^{2\si+2}\lesssim 1
  + y(t) 
  + \sup_{0\le s\le t}y(s)\lesssim 1 +\sup_{0\le s\le t}y(s).  
\end{equation*}
Gronwall lemma implies $f\in
L^\infty([0,M])$. We infer $y\in L^\infty_{\rm loc}(\R)$. With the above
inequality, this implies $\nabla u \in L^\infty_{\rm loc}(\R;L^2)$, and
Theorem~\ref{theo:global} then follows from
Corollary~\ref{cor:alternative}.

\section{Generalized Mehler formula}
\label{sec:mehler}

In the rest of this paper, we consider the case where $V$ is exactly
quadratic in $x$, and study some properties associated to
\eqref{eq:NLSharmo}. 

\subsection{The formula}
\label{sec:2.1}

Classically, Mehler's formula refers to the explicit formula for the
fundamental solution of the linear equation 
\begin{equation}
  \label{eq:linear}
   i\d_t u_{\rm lin} + \frac{1}{2}\Delta u_{\rm lin} =
   \frac{1}{2}\sum_{j=1}^d\Om_j(t)x_j^2 
  u_{\rm lin} \quad ;\quad u_{{\rm lin}\mid t=0}=u_0,
\end{equation}
in the case $\dot\Om_j=0$, with $\Om_j>0$. See
e.g. \cite{Feyn}. It was generalized (still with $\dot\Om_j=0$)
in  \cite{HormanderQuad} to a framework where typically, 
$\Om_j\in \R$ has no specified sign. 
\smallbreak

The case of time dependent harmonic potentials with $d=1$ was
considered in \cite{CLSS08}, along with other terms corresponding for
instance to  time dependent magnetic and electric fields. Since the
case $d\ge 1$ for \eqref{eq:linear} follows by taking the tensor
product of the one dimensional case, we shall simply rewrite the
results of  \cite{CLSS08} (and adapt them to our conventions). 
\smallbreak

Seek formally the solution to \eqref{eq:linear} as 
  \begin{equation}\label{eq:mehler}
  u_{\rm lin}(t,x)= \(\prod_{j=1}^d\frac{1}{ 2i\pi \mu_j(t)}\)^{1/2}\int_{\R^d}
  e^{\frac{i}{2}\phi(t,x,y)}u_0(y)dy,
\end{equation}
where
\begin{equation*}
  \phi(t,x,y) =\sum_{j=1}^d \(\alpha_j(t)x_j^2 + 2\beta_j(t)x_jy_j +
  \gamma_j(t) y_j^2 +2\delta_j(t)x_j +2\epsilon_j(t)y_j\) +\theta(t), 
\end{equation*}
and all the functions of time involved in this formula are
real-valued. For instance, when $\Om=0$, we have $\mu(t)=t$, $\alpha =
\beta= \gamma=1/t$ and $\delta=\epsilon=\theta=0$: the convergence
$u_{\rm lin}(t)\to 
u_0$ as $t\to 0$ is recovered (at least formally) by applying the
stationary phase formula. Note that in view of the results of
D.~Fujiwara \cite{Fujiwara79,Fujiwara}, we know that there exists
$\eta>0$ such that for $|t|<\eta$, the solution to
\eqref{eq:linear} can be expressed as
\begin{equation*}
  u_{\rm lin}(t,x) = \frac{1}{(2i\pi t)^{d/2}}\int_{\R^d}
  e^{i\varphi(t,x,y)}a(t,x,y)u_0(y)dy, 
\end{equation*}
where $a(0,x,y)=1$, $\d_x^\alpha\d_y^\beta a\in
L^\infty(]-\eta,\eta[\times\R^d\times \R^d)$ for all $\alpha,\beta\in
\N^d$, and
\begin{equation*}
  \varphi(t,x,y) = \frac{|x-y|^2}{2t} + t\xi(t,x,y),
\end{equation*}
with $\d_x^\alpha\d_y^\beta \xi\in
L^\infty(]-\eta,\eta[\times\R^d\times \R^d)$ as soon as
$|\alpha+\beta|\ge 2$.

\smallbreak

Applying the differential operator $i\d_t
+\frac{1}{2}\Delta$ to \eqref{eq:mehler}, and identifying the terms (in
$x_j^2,x_jy_j\ldots$) in \eqref{eq:linear}, we find:
\begin{align*}
  x_j^2:& \quad \dot \alpha_j+\alpha_j^2+\Om_j=0;&& x_jy_j:\quad
  \dot\beta_j+\alpha_j\beta_j=0.\\
y_j^2:&\quad \dot \gamma_j + \beta_j^2=0;&& x_j:\quad \dot \delta_j +
\alpha_j\delta_j = 0.\\
y_j:&\quad \dot\epsilon_j +\beta_j \delta_j =0;&& \IM(\C): \quad \dot
\mu_j= \alpha_j \mu_j.\\ 
\RE(\C):&\quad \dot
\theta+\sum_{j=1}^d\delta_j^2 =0. 
\end{align*}
We infer that $\mu_j$ is given by
\begin{equation}
  \label{eq:mu}
  \ddot \mu_j +\Om_j(t)\mu_j =0\quad ;\quad \mu_j(0)=0,\quad \dot\mu_j(0)=1.
\end{equation}
We also have
\begin{equation*}
  \alpha_j
  =\frac{\dot \mu_j}{\mu_j}. 
\end{equation*}
Note that as in the standard cases ($\dot \Om_j=0$), 
$\alpha_j(t)\sim 1/t$ as $t\to 0$. For $\beta_j$, we have
\begin{equation*}
  \dot \beta_j+ \frac{\dot \mu_j}{\mu_j}
  \beta_j=0,\text{ hence }
  \beta_j(t)=\frac{C}{\mu_j(t)}, 
\end{equation*}
and the stationary phase formula (as $t\to 0$) yields  $C=-1$. 
We also find
\begin{equation*}
  \gamma_j (t)= \frac{1}{\mu_j(t)\dot\mu_j(t)} -\int_0^t
  \frac{\Om_j(\tau)}{\(\dot     \mu_j(\tau)\)^2}d\tau.
\end{equation*}
Since $\delta_j(0)=\epsilon_j(0)=\theta(0)=0$, we have
$\delta_j=\epsilon_j=\theta_j\equiv 0$. 
\begin{remark}
  The case of the usual harmonic potential ($\Om_j=1$) shows that
  singularities may be present in the fundamental solutions for
  positive times, corresponding to the zeroes of $\mu_j$; see
  e.g. \cite{CKS95,KRY,Yajima96,Zelditch83}.  
\end{remark}
\begin{remark}
  The dispersive properties associated to \eqref{eq:linear} are
  measured by the $\mu_j$'s. We will see for instance that if
  $\Om_j\le 0$ for all $j$, then global in time Strichartz estimates
  are available, as in the case $\Om_j=0$.   
\end{remark}
To summarize, we have:
\begin{lemma}\label{lem:mehler}
  Let $d\ge 1$, and $\Om_j \in C(\R;\R)$ be locally
  Lipschitzean. There exists $T>0$ such 
  that for $u_0\in \Sch(\R^d)$,
  the solution to \eqref{eq:linear} is given, for $|t|<T$, by:
 \begin{equation*}
  u_{\rm lin}(t,x)= \(\prod_{j=1}^d\frac{1}{ 2i\pi \mu_j(t)}\)^{1/2}\int_{\R^d}
  e^{\frac{i}{2}\sum_{j=1}^d \(\alpha_j(t)x_j^2 + 2\beta_j(t)x_jy_j +
  \gamma_j(t) y_j^2 \)}u_0(y)dy,
\end{equation*}
where
\begin{align*}
  &\ddot \mu_j +\Om_j(t)\mu_j =0\quad ;\quad \mu_j(0)=0,\quad
  \dot\mu_j(0)=1,\\
& \alpha_j =\frac{\dot \mu_j}{\mu_j}\quad ;\quad \beta_j =
-\frac{1}{\mu_j}\quad ;\quad \gamma_j (t)=
\frac{1}{\mu_j(t)\dot\mu_j(t)} -\int_0^t 
  \frac{\Om_j(\tau)}{\(\dot \mu_j(\tau)\)^2}d\tau.
\end{align*}
\end{lemma}

\begin{remark}
  The fact that the quadratic potential has no rectangle term is not
  necessary in order to get such a result. If we consider
\begin{equation*}
   i\d_t u_{\rm lin} + \frac{1}{2}\Delta u_{\rm lin} =
   \frac{1}{2}\<M(t)x,x\>u_{\rm lin} ,
\end{equation*}
where $M(t)$ is a (time dependent) symmetric matrix, then a similar
formula is available. Of course, the formula is more involved, and since it
does not really bring new information, we do not carry out the
computation here. 
\end{remark}

\subsection{Some consequences}
\label{sec:2.2}

In this paragraph, we assume that the functions $\Om_j$ are
bounded. This assumption was discussed in
Example~\ref{ex:StrichartzUSTC}.

As a consequence of the boundedness of $\Om_j$, we infer a uniform
local bound from below for the functions $\mu_j$. It follows from
the growth of the functions $\mu_j$'s, which is at most exponential:
\begin{lemma}\label{lem:munuexp}
  Assume that for all $j\in\{1,\ldots,d\}$, $\Om_j\in C(\R;\R)$
  is locally  Lipschitzean and \emph{bounded}. For $s\in \R$, define
  $\mu_j^s$ and $\nu_j^s$ as the  solutions to 
\begin{align}
\label{eq:mus}
 & \ddot \mu_j^s+\Om_j(t) \mu_j^s=0\quad;\quad \mu_j^s(s)=0,\quad \dot
  \mu_j^s(s)=1.\\
  & \ddot \nu_j^s+\Om_j(t) \nu_j^s=0\quad;\quad \nu_j^s(s)=1,\quad \dot
  \nu_j^s(s)=0.\label{eq:nus}
\end{align}
There exists $C>0$ independent of $s\in \R$ such that 
\begin{equation*}
  |\mu_j^s(t)|+|\dot\mu_j^s(t)| + |\nu_j^s(t)|+|\dot\nu_j^s(t)|\le C
  e^{C|t-s|},\quad \forall t\in \R. 
\end{equation*}
\end{lemma}
\begin{proof}
Introduce $f_j^s(t) = |\dot\mu_j^s(t)|+|\mu_j^s(t)|$. We have
\begin{align*}
  \dot f_j^s(t) &\le |\ddot \mu_j^s(t)|+|\dot\mu_j^s(t)| = \lvert \Om_j(t)
  \mu_j^s(t)\rvert + |\dot\mu_j^s(t)|\\
& \le
  \|\Om_j\|_{L^\infty}|\mu_j^s(t)| + |\dot\mu_j^s(t)| \lesssim f_j^s(t).  
\end{align*}
  Gronwall lemma yields, since $f_j^s(s)=1$,
\begin{equation*}
 f_j^s(t) \lesssim e^{C|t-s|},
\end{equation*}
for some $C>0$,   independent of $j,s$ and $t$.
The first part of lemma then follows. The second estimate is
similar. 
\end{proof}
In view of the initial data for $\mu_j^s$ and $\nu_j^s$, we infer:
\begin{lemma}\label{lem:dispUSTC}
  Assume that for all $j\in\{1,\ldots,d\}$, $\Om_j\in C(\R;\R)$
  is locally Lipschitzean and \emph{bounded}. There exists $\eta>0$
  such that for all $j$, and all  $s\in \R$, 
  \begin{equation*}
    |\mu_j^s(t)|\ge \frac{|t-s|}{2},\quad
    \frac{1}{2}\le|\nu^s_j(t)|\le \frac{3}{2},\quad \forall t,\
    |t-s|<\eta, 
  \end{equation*}
where $\mu_j^s$ and $\nu^s_j$ are given by \eqref{eq:mus} and
\eqref{eq:nus}, respectively. 
\end{lemma}
This yields a uniform local dispersion in \eqref{eq:disp}, and 
we infer a property which will be crucial in the study of the large
time behavior of high Sobolev norms:
\begin{proposition}\label{prop:strichartz2}
  Assume that
for all $j$, $\Om_j\in C(\R;\R)$
  is locally Lipschitzean and \emph{bounded}. Then
  Proposition~\ref{prop:strichartz} remains   valid with $T=\infty$. 
\end{proposition}

\section{Generalized lens transform}
\label{sec:lens}
\subsection{The formula}
\label{sec:formula}
It was noticed in \cite{KavianWeissler} that in the case of the
$L^2$-critical nonlinearity ($\si=2/d$), an explicit change of unknown
function makes it possible to add or remove an \emph{isotropic}
harmonic potential: if $v$ solves
\begin{equation}
  \label{eq:NLSconf}
  i\d_t v+\frac{1}{2}\Delta v = \lambda \lvert v\rvert^{4/d}v\quad
  ;\quad v_{\mid t=0}=u_0,
\end{equation}
where $\lambda \in \R$, then $u$, given for $|t|<\pi/(2\om)$ by the
lens transform
\begin{equation}\label{eq:lens}
  u(t,x) = \frac{1}{\(\cos(\om t)\)^{d/2}}v\(\frac{\tan(\om
    t)}{\om},\frac{x}{\cos(\om t)}\)e^{-i\frac{\om }{2}|x|^2\tan (\om t)} 
\end{equation}
solves
\begin{equation*}
  i\d_t u + \frac{1}{2}\Delta u = \frac{\om^2}{2}\lvert x\rvert^2 u +
  \lambda \lvert u\rvert^{4/d}u\quad ;\quad u_{\mid t=0}=u_0. 
\end{equation*}
See also \cite{Rybin,CaM3AS,TaoLens}. Note that the change for the
time variable is locally invertible, not globally. The case of a
repulsive harmonic potential,
\begin{equation*}
  i\d_t u + \frac{1}{2}\Delta u = -\frac{\om^2}{2}\lvert x\rvert^2 u +
  \lambda \lvert u\rvert^{4/d}u\quad ;\quad u_{\mid t=0}=u_0,
\end{equation*}
is obtained by replacing $\om$ by $i\om$: a formula similar to
\eqref{eq:lens} follows, where the trigonometric functions are
replaced by hyperbolic functions (and the discussion on the time
interval becomes different), see \cite{CaSIMA}. A heuristic way to
understand why this approach works only in the case of isotropic
potentials is that even though there would be a ``natural'' candidate
to change the space variable in the anisotropic case, there is no
satisfactory candidate to change the time variable. 
\smallbreak

The lens transform can be generalized to the case of
\eqref{eq:NLSharmo} provided that the potential is isotropic in the
sense that $\Om_j(t)=\Om(t)$ is independent of $j$. Seek an
extension of \eqref{eq:lens} of the form
\begin{equation}
  \label{eq:lensgen}
  u(t,x)=\frac{1}{b(t)^{d/2}}v\(\zeta(t),\frac{x}{b(t)}\)
e^{\frac{i}{2}a(t)|x|^2}, 
\end{equation}
with $a,b,\zeta$ real-valued, 
\begin{equation}\label{eq:cilens}
  b(0)=1\quad ;\quad 
a(0)=\zeta(0)=0. 
\end{equation}
Suppose also that $v$ solves a more general
non-autonomous equation
\begin{equation}
  \label{eq:nlsnonauto}
  i\d_t v+\frac{1}{2}\Delta v = H(t)\lvert v\rvert^{2\si}v\quad
  ;\quad v_{\mid t=0}=u_0.
\end{equation}
We want $u$ to solve
\begin{equation}
  \label{eq:unonauto}
  i\d_t u + \frac{1}{2}\Delta u = \frac{1}{2}\Om(t)\lvert x\rvert^2
  u + h(t) \lvert u\rvert^{2\si}u\quad ;\quad u_{\mid t=0}=u_0.
\end{equation}
Apply the Schr\"odinger differential operator to the formula
\eqref{eq:lensgen}, and identify the terms so that $u$ solves
\eqref{eq:unonauto}. We find: 
\begin{equation*}
  \dot b = ab\quad ;\quad \dot a+a^2 +\Om=0\quad ;\quad
\dot \zeta = \frac{1}{b^2}\quad 
;\quad b(t)^{d\si-2}H\(\zeta(t)\)=h(t). 
\end{equation*}
Introduce the solution to 
\begin{equation}
  \label{eq:solfond}
\left\{
  \begin{aligned}
&\ddot \mu + \Om(t)\mu =0 \quad;\quad \mu(0)=0,\quad
  \dot\mu(0)=1. \\
 &   \ddot \nu + \Om(t)\nu =0 \quad;\quad \nu(0)=1,\quad
  \dot\nu(0)=0.  
  \end{aligned}
  \right.
\end{equation}
Note that since the Wronskian of $\mu$ and $\nu$ is constant, we have
$\dot\mu\nu- \mu\dot\nu=1$ for all time. This relation extends the
identities $\cos^2 t+\sin^2 t=1$ and
$\cosh^2t-\sinh^2t=1$. 
In view of \eqref{eq:cilens}, we infer:
\begin{equation*}
  a=\frac{\dot \nu}{\nu}\quad ;\quad b= \nu\quad ;\quad
  \zeta=\frac{\mu}{\nu}. 
\end{equation*}
Note that $\zeta$ is locally invertible, since $\zeta(0)=0$ and 
\begin{equation*}
  \dot \zeta = \frac{1}{b^2}=\frac{1}{\nu^2},\text{ hence }\dot \zeta(0)=1.
\end{equation*}
Therefore, the lens transform is locally
invertible. Moreover, since $b(0)=\nu(0)=1$, we can write, locally in
time,
\begin{equation*}
  H(t) = b\(\zeta^{-1}(t)\)^{2-d\si}h\(\zeta^{-1}(t)\)=
  \nu\(\(\frac{\mu}{\nu}\)^{-1}(t)\)^{2-d\si}h\(\(\frac{\mu}{\nu}\)^{-1}(t)\). 
\end{equation*}
\begin{proposition}\label{prop:lens}
  Let $v$ solve
\begin{equation*}
  i\d_t v+\frac{1}{2}\Delta v = H(t)\lvert v\rvert^{2\si}v\quad
  ;\quad v_{\mid t=0}=u_0.
\end{equation*}
Let $\Om\in C(\R;\R)$. There exists $T>0$ such that the
following holds. Define $u$ by
 \begin{equation*}
  u(t,x)=\frac{1}{\nu(t)^{d/2}}v\(\frac{\mu(t)}{\nu(t)},\frac{x}{\nu(t)}\)
e^{i\frac{\dot \nu(t)}{\nu(t)}\frac{|x|^2}{2}},\quad |t|\le T,
\end{equation*}
where $(\mu,\nu)$ is given by \eqref{eq:solfond}. 
Then for $|t|<\mu(T)/\nu(T)$, $u$ solves 
\begin{equation*}
 i\d_t u + \frac{1}{2}\Delta u = \frac{1}{2}\Om(t)\lvert x\rvert^2
  u + h(t) \lvert u\rvert^{2\si}u\quad ;\quad u_{\mid t=0}=u_0,
\end{equation*}
where $h(t)= \nu(t)^{d\si-2}H\(\mu(t)/\nu(t)\)$. 
\end{proposition}
\begin{remark}
  We do not require $\Om$ to be locally Lipschitzean: all we need is
  the local existence of a $C^2$ solution to \eqref{eq:solfond}, so we can
  rely on Peano existence theorem. 
\end{remark}
\begin{remark}[Generalized Avron--Herbst formula]
  A similar result is available, when the quadratic potential
  $\Om(t)|x|^2$ is replaced by a linear (anisotropic) one $E(t)\cdot
  x$, where $E\in C(\R;\R^d)$. The solutions to 
  \begin{align*}
    &i\d_t v+\frac{1}{2}\Delta v = h(t) |v|^{2\si}v\quad ;\quad v_{\mid
      t=0}=u_0,\\
&i\d_t u +\frac{1}{2}\Delta u = E(t)\cdot x u +h(t) |u|^{2\si}u \quad
;\quad u_{\mid       t=0}=u_0, 
  \end{align*}
are related by the formula
\begin{equation*}
  u(t,x)=v\(t,x+\int_0^t \int_0^\tau E(s)ds d\tau\) e^{-ix\cdot
    \int_0^tE(\tau)d\tau  -\frac{i}{2}\int_0^t|E(\tau)|^2d\tau}.
\end{equation*}
\end{remark}
\subsection{Proof of Proposition~\ref{prop:blowup}}
 We assume in this paragraph that the
nonlinearity is \emph{focusing}: $\l<0$. By homogeneity, we can assume
$\l=-1$. It is well known that the equation
\begin{equation}\label{eq:nlsfoc}
  i\d_t v+\frac{1}{2}\Delta v=-|v|^{4/d}v
\end{equation}
possesses solutions which blow up in finite time, with different
possible rates (see
e.g. \cite{BourgainLivre,CazCourant,Rap06,Sulem,TaoDisp} and
references therein).  
\smallbreak

By adapting Proposition~\ref{prop:lens} to isolate the initial time
$t=0$, we see that the lens transform maps a solution to
\eqref{eq:nlsfoc} which blows up at time $t=0$ to a solution to
\begin{equation}\label{eq:nlsharmofoc}
  i\d_t u +\frac{1}{2}\Delta u = \frac{1}{2}\Om(t)|x|^2 u -|u|^{4/d}u
\end{equation}
which blows up at time $t=0$. Note that the blow-up rate is not
altered by the lens transform, since $\nu(t)\approx 1$ and
$\mu(t)/\nu(t)\approx t$ as $t\to 0$. 
\smallbreak

Typically, consider the (unstable) minimal mass blow-up solution to
\eqref{eq:nlsfoc}:
\begin{equation*}
  v(t,x) =\frac{1}{t^{d/2}}Q\(\frac{x}{t}\) e^{i\frac{|x|^2}{2t} -\frac{i}{t}},
\end{equation*}
where $Q$ is the ground state, defined as the unique positive, radial,
solution to
\begin{equation*}
  -\frac{1}{2}\Delta Q+Q = Q^{1+4/d}. 
\end{equation*}
 The lens transform yields a corresponding blow-up solution to
 \eqref{eq:nlsharmofoc} given by 
\begin{equation*}
  u(t,x)= \frac{1}{\mu(t)^{d/2}}Q\(\frac{x}{\mu(t)}\)e^{i\frac{\dot
      \mu(t)}{\mu(t)}\frac{|x|^2}{2} -i\frac{\nu(t)}{\mu(t)}}.
\end{equation*}
To our knowledge, this gives the first example of an explicit blow-up
solution in the presence of a time-dependent external potential.

Note that we have considered the explicit case of minimal mass blow-up
solutions for convenience. Any blow-up solution for \eqref{eq:nlsfoc}
gives rise to a blow-up solution for \eqref{eq:nlsharmofoc},
\emph{with the same blow-up rate}. 

Note also that without extra assumption on $\Om$, the Sobolev norms of
$u$ may have an arbitrary growth as $t\to \infty$. 

\begin{example}
 Consider  $\mu(t) = \exp\(1-e^t\)-\exp\(1-e^{2t}\)$ (which satisfies
 $\mu(0)=0$ and
  $\dot \mu(0)=1$). Then the growth of Sobolev norms of the function
  $u$ given by the above formula is given by a double exponential in
  time, since
\begin{equation*}
  \|u(t)\|_{H^s}\Eq t{+\infty} \frac{C_s}{|\mu(t)|^s}.
\end{equation*}
To determine the corresponding function $\Om$, we  compute
  \begin{equation*}
    \ddot \mu(t)=
  \(e^{2t}-e^t\)\(\exp\(1-e^t\)-4\exp\(1-e^{2t}\)\),
  \end{equation*}
and therefore
  \begin{equation*}
    \Om(t) =
    \frac{\exp\(1-e^t\)-4\exp\(1-e^{2t}\)}{\exp\(1-e^t\)-\exp\(1-e^{2t}\)}
\(e^t-e^{2t}\).  
  \end{equation*}
We note that $\Om(t)\sim -e^{2t}$ as $t\to +\infty$: the harmonic
potential is repulsive ($\Om<0$), and becomes exponentially stronger as time
increases. 
\end{example}

\section{Vector fields}
\label{sec:champs}
The aim of this paragraph is to show that there exists vector fields
which may be useful to study the nonlinear equation
\eqref{eq:NLSharmo}, in the same spirit as in \cite{CaSIMA,CaCCM}. 
Consider the solutions to
\begin{equation}
  \label{eq:solfondj}
\left\{
  \begin{aligned}
&\ddot \mu_j + \Om_j(t)\mu_j =0 \quad;\quad \mu_j(0)=0,\quad
  \dot\mu_j(0)=1. \\
 &   \ddot \nu_j + \Om_j(t)\nu_j =0 \quad;\quad \nu_j(0)=1,\quad
  \dot\nu_j(0)=0. 
  \end{aligned}
  \right.
\end{equation}
We define 
\begin{align*}
  A_j &= \dot \mu_j x_j +i\mu_j \d_j = i\mu_j
  e^{i\frac{x_j^2}{2}\frac{\dot
      \mu_j}{\mu_j}}\d_j\(e^{-i\frac{x_j^2}{2}\frac{\dot \mu_j}{\mu_j}}
  \cdot\)=i\mu_j
  e^{i\sum_k\frac{x_k^2}{2}\frac{\dot
      \mu_k}{\mu_k}}\d_j\(e^{-i\sum_k\frac{x_k^2}{2}\frac{\dot \mu_k}{\mu_k}}
  \cdot\),\\
B_j &= \dot \nu_j x_j +i\nu_j \d_j = i\nu_j
  e^{i\frac{x_j^2}{2}\frac{\dot
      \nu_j}{\nu_j}}\d_j\(e^{-i\frac{x_j^2}{2}\frac{\dot \nu_j}{\nu_j}}
  \cdot\)= i\nu_j
  e^{i\sum_k\frac{x_k^2}{2}\frac{\dot
      \nu_k}{\nu_k}}\d_j\(e^{-i\sum_k\frac{x_k^2}{2}\frac{\dot \nu_k}{\nu_k}}
  \cdot\).
\end{align*}
Note that the last two expressions for $A$ or $B$ show that $A$ and
$B$ act on gauge invariant nonlinearities like derivatives: the
modulus ignores the multiplication by the exponential. 
\begin{example}[$\dot \Om_j=0$]
  When $\Om_j=0$, $B_j=i\d_j$ and $A_j = x_j+it\d_j$, which are
  Heisenberg derivatives commonly used in the theory of nonlinear
  Schr\"odinger equations (see e.g. \cite{CazCourant}). When $\Om_j=\om_j^2>0$,
  $A_j = x_j\cos (\om_j t) +i\frac{\sin (\om_j t)}{\om_j} \d_j$ and
  $B_j= \om_jx_j\sin(\om_j t)+i\cos(\om_j t)\d_j$: we recover
  classical Heisenberg derivatives (see e.g. \cite{Cycon}). In
  these two cases (as well as in the case $\Om_j=-\om_j^2<0$), we have
  \begin{equation*}
    A_j = U_V(t)x_j U_V(-t)\quad ;\quad B_j = U_V(t)i\d_j U_V(-t),
  \end{equation*}
where $U_V(t)=\exp\(-it(-\frac{1}{2}\Delta+V(x))\)$, $V(x)=\sum_{k=1}^d
\Om_k x_k^2$. 
\end{example}
More generally, consider $\dot\eta_j x_j+i\eta_j \d_j$: we
check that this 
operator commute with the linear operator
\begin{equation*}
  i\d_t +\frac{1}{2}\Delta -\frac{1}{2}\sum_{k=1}^d \Om_k(t)x_k^2
\end{equation*}
if and only if $\eta_j$  satisfies $\ddot \eta_j+ \Om_j \eta_j=0$:
\begin{align*}
  \left[i\d_t +\frac{1}{2}\Delta -\frac{1}{2}\sum_{k=1}^d
    \Om_k(t)x_k^2, \eta_j x_j+i\eta_j \d_j\right]
&= \left[i\d_t +\frac{1}{2}\d_j^2-\frac{1}{2}
    \Om_j(t)x_j^2, \eta_j x_j+i\eta_j \d_j\right]\\
&= i\ddot \eta_j x_j -\dot \eta_j
\d_j +\dot \eta_j\d_j +i\eta_j\Om_j x_j. 
\end{align*}
This is zero if (and only if) $\ddot \eta_j+ \Om_j \eta_j=0$.
\begin{remark}
  This computation could be extended to the case where the center of
  the harmonic potential depends on time:
  \begin{equation*}
    i\d_t u+ \frac{1}{2}\Delta u =
    \frac{1}{2}\sum_{k=1}^d\Om_k(t)\(x_k-c_k(t)\)^2 u .
  \end{equation*}
Replacing $\dot\eta_j x_j+i\eta_j \d_j$ with $\dot\eta_j
\(x_j-y_j(t)\)+i\eta_j \d_j$, we can repeat the above computation, and
check that the two operators commute if and only if
$\ddot \eta_j +\Om_j \eta_j=0$ and $\ddot \eta_j y_j + \dot \eta_j
\dot y_j +\eta_j \Om_j c_j=0$. We choose not to investigate this case
further into details here. 
\end{remark}
To show that the $\Sigma$-norm of $u$ is related to the $L^2$-norms of
$A_ju$ and $B_ju$, write
\begin{equation*}
  \( 
  \begin{array}[c]{c}
    A_j\\
B_j
  \end{array}
\) = M_j   \( 
  \begin{array}[c]{c}
    x_j\\
i\d_j
  \end{array}
\) ,\quad \text{where }M_j =
\(
\begin{array}[c]{cc}
  \dot\mu_j & \mu_j\\
\dot\nu_j & \nu_j
\end{array}
\).
\end{equation*}
We note that the determinant of $M_j$ is the Wronskian of $\mu_j$ and
$\nu_j$:
\begin{equation*}
  \operatorname{det}M_j = \nu_j\dot\mu_j-\mu_j\dot\nu_j\equiv 1.
\end{equation*}
Therefore
\begin{equation}\label{eq:inv}
   \( 
  \begin{array}[c]{c}
    x_j\\
i\d_j
  \end{array}
\) = \(
\begin{array}[c]{cc}
  \nu_j & -\mu_j\\
-\dot\nu_j & \dot\mu_j
\end{array}
\)   
\( 
  \begin{array}[c]{c}
    A_j\\
B_j
  \end{array}
\). 
\end{equation}
We shall use these vector fields in the isotropic case, where they
provide \emph{a priori} estimates:
\begin{equation}\label{eq:isotrop}
  i\d_t u+\frac{1}{2}\Delta u = \frac{1}{2}\Om(t)|x|^2 u +\l
  |u|^{2\si}u. 
\end{equation}
Since $A$ commutes with the linear part of \eqref{eq:isotrop} and acts on
gauge invariant nonlinearities like a gradient, we have readily
\begin{equation*}
 \frac{1}{2} \frac{d}{dt}\|Au\|_{L^2}^2 =
 \l\si\IM\int_{\R^d}|u|^{2\si-2}u^2 \(\overline{ A u}\)^2. 
\end{equation*}
Expanding $(Au)^2$, we obtain eventually:
\begin{align}
  \frac{d}{dt}\(\frac{1}{2}\|Au\|_{L^2}^2 + \frac{\l \mu^2}{\si+1}
  \|u\|_{L^{2\si+2}}^{2\si+2}\)&= \frac{\l}{\si+1}\mu \dot
  \mu\(2-d\si\) \|u\|_{L^{2\si+2}}^{2\si+2},
\label{eq:A}\\
\frac{d}{dt}\(\frac{1}{2}\|Bu\|_{L^2}^2 + \frac{\l \nu^2}{\si+1}
  \|u\|_{L^{2\si+2}}^{2\si+2}\)&= \frac{\l}{\si+1}\nu \dot
  \nu\(2-d\si\) \|u\|_{L^{2\si+2}}^{2\si+2}.\label{eq:B}
\end{align}
These evolution laws are the analogue of the pseudo-conformal
conservation law (see \cite{CaCCM} for the case $\dot \Om=0$). They
will allow us to infer scattering results in
the case $\Om\le 0$, $\l\ge 0$ (\S\ref{sec:isorep}).

\section{Growth of higher order Sobolev norms and momenta}
\label{sec:growth}

\subsection{The linear case}
\label{sec:growthlinear}

In this paragraph, we assume $\l=0$. We recall that in general,
Mehler's formula is valid only locally in 
time, since singularities may appear in the fundamental solution; see
e.g. \cite{CKS95,KRY,Yajima96,Zelditch83}. To understand the long time
behavior of the solution $u_{\rm lin}$ to \eqref{eq:linear}, one may
use Egorov Theorem (see e.g. \cite{BR02}). Since we deal with a
time-dependent potential, modifications would be needed in Egorov
Theorem, and we rather follow another strategy, to have some estimates
in the linear case (instead of an exact asymptotic behavior, as Egorov
Theorem would give us). This approach is based on the vector fields
introduced in \S\ref{sec:champs}.  
\smallbreak

We remark that since the $L^2$-norm of $u_{\rm lin}$ does not depend
on time, and since the operator $A_j$ and $B_j$ introduced in
\S\ref{sec:champs} commute with Equation~\eqref{eq:linear}, the
$L^2$-norm of $A_{j_1}B_{j_2}\ldots A_{j_k}u_{\rm lin}$ is constant,
for whichever combination of these vector fields. In view of
\eqref{eq:inv}, we infer, for $k\in\N$,
\begin{equation*}
  \lVert \lvert x\rvert^k u_{\rm lin}(t)\rVert_{L^2} + \lVert u_{\rm
    lin}(t)\rVert_{H^k} \lesssim \sum_{j=1}^d\(\lvert \mu_j (t)\rvert^k + 
  \lvert \nu_j (t)\rvert^k \).
\end{equation*}
Lemma~\ref{lem:munuexp} shows that if $\Om_j \in C(\R;\R)$ is locally
Lipschitzean and bounded, then the above quantity grows at most
exponentially in time. By Proposition~\ref{prop:strichartz2}, we
conclude that $(Exp)_k$ is satisfied for all $k$, provided $u_0$ is
sufficiently smooth and localized. We recall that the case $\Om_j=-1$
shows that the exponential growth may occur, and that in $(Exp)_k$,
the constant $C$ must be expected to depend on $k$ ($C=k$ when
$\Om=-1$ is sharp).

\subsection{The $L^2$-subcritical case}
\label{sec:growthgen}
\begin{lemma}\label{lem:growthgeneral}
  Let $\sigma,k\in \N$, with $\si\le 2/d$, $\Om_j\in
  C(\R;\R)$ be locally Lipschitzean and \emph{bounded},  and  
  \begin{equation*}
    u_0 \in H^k\(\R^d\),\text{ with }|x|^ku_0\in L^2(\R^d). 
  \end{equation*}
Suppose that there exists $f\in C(\R_+;\R_+)$ with $f(0)=0$ such that 
\begin{equation}
  \label{eq:unifLp}
  \|u\|_{L^\theta([s,s+\tau];L^q)}\le f(\tau),\quad \forall s,\tau\in \R,
\end{equation}
where 
\begin{equation*}
  q=2\si+2\quad ;\quad
  \theta=\frac{2\si(2\si+2)}{2-(d-2)\si}.
\end{equation*}
Then the solution to \eqref{eq:NLSharmo} satisfies $(Exp)_k$.
\end{lemma}
\begin{proof}
  The first step consists in resuming the computations carried out in
  the proof of Proposition~\ref{prop:Sigmaloc2}, in the case
  $k=0$. The case $k\ge 1$ will follow by induction (recall that
  the constant $C$ in the exponential growth must be expected to
  depend on $k$). \\
{\bf Case $k=0$.} Let us pretend that the $L^2$-norm of $u$ is not
conserved, to simplify the induction. Resuming the same numerology as
in the proof of 
Proposition~\ref{prop:Sigmaloc2}, Strichartz estimates yield, for all
$t\in \R$ and $\tau>0$,
\begin{align*}
  \|u\|_{L^p([t,t+\tau];L^q)\cap L^\infty([t,t+\tau];L^2)}&\lesssim
  \|u(t)\|_{L^2} + \|u\|_{L^\theta([t,t+\tau];L^q)}^{2\si}
  \|u\|_{L^p([t,t+\tau];L^q)}\\
& \lesssim
  \|u(t)\|_{L^2} + f(\tau)^{2\si}
  \|u\|_{L^p([t,t+\tau];L^q)}.
\end{align*}
Fix $\tau\ll 1$ once and for all so the last term of the right hand
side can be absorbed by the left hand side, up to doubling the
estimating constant: at every increment of time of length $\tau$, the
$L^2$ norm is multiplied (at most) by some fixed constant $C$. This
implies that it grows at most exponentially. Using Strichartz
estimates again, we conclude that $(Exp)_0$ is satisfied (and actually,
$(Alg)_0$ is also true). 
\smallbreak

\noindent {\bf Case $k\ge 1$.} For $k\ge 1$, suppose that
$(Exp)_{k-1}$ is satisfied. To avoid a lengthy presentation, we denote
by $w_\ell$ the family of combinations of $\alpha$ momenta and $\beta$
space derivatives of $u$, with $|\alpha|+|\beta|=\ell$ ($w_0=u$). We have,
rather formally,
\begin{equation}\label{eq:wk}
  i\d_t w_k +\frac{1}{2}\Delta w_k = \frac{1}{2}\sum_{j=1}^d
  \Om_j(t)x_j^2 w_k + {\tt V}(u,w_k) + F +L(w_k),
\end{equation}
where ${\tt V}$ is homogeneous of degree $2\si$ with respect to its
first argument, 
$\R$-linear with respect to its second argument, $F$ satisfies the
pointwise estimate
\begin{equation*}
  |F|\lesssim \sum_{0\le \ell_j\le k-1}\lvert w_{\ell_1}\rvert\ldots
  \lvert w_{\ell_{2\si+1}}\rvert, 
\end{equation*}
where the sums carries over combinations such that in addition $\sum
\ell_j=k$ ($F=0$ in the case $k=1$), and $L$ is linear with respect to
its argument. A word of 
explanation is needed about $L$: this term stems from the fact that
$x$ and $\nabla$ do not commute with the linear part of the equation. One
might argue that we could proceed as in the linear case, and use
repeatedly the vector fields $A_j$ and $B_j$. The problem is that even
though $A_j$ and $B_j$ act on gauge invariant nonlinearities like
derivatives, this is not so, for instance, for $A_jB_j$ (the phases do
not cancel in the factored formula). We might use the operators
$A_{j_1}\ldots A_{j_k}$ and $B_{j_1}\ldots B_{j_k}$, but this does not
suffice to recover the momenta and derivatives of $u$, since
``rectangle'' terms (like $A_jB_j$) would be needed. 
\smallbreak

We proceed in the same spirit as in the case $k=0$:
\begin{align*}
  \|w_k\|_{L^p([t,t+\tau];L^q)\cap L^\infty([t,t+\tau];L^2)}&\lesssim
  \|w_k(t)\|_{L^2} + \|u\|_{L^\theta([t,t+\tau];L^q)}^{2\si}
  \|w_k\|_{L^p([t,t+\tau];L^q)}\\
+ \sum_{0\le \ell_j\le
    k-1}\|w_{\ell_1}\|_{L^\theta([t,t+\tau];L^q)}&\ldots
  \|w_{\ell_{2\si}}\|_{L^\theta([t,t+\tau];L^q)}
\|w_{\ell_{2\si+1}}\|_{L^p([t,t+\tau];L^q)} 
  \\ 
& \quad + \|L(w_k)\|_{L^1([t,t+\tau];L^2)}
\end{align*}
Fixing $\tau\ll 1$ independent of $t\in \R$, the second term of the
right hand side is absorbed by the left hand side. The sum is treated
thanks to $(Exp)_{k-1}$. We notice that since $\si\le 2/d$, we have
$\theta\le p$, where we recall that $(p,q)$ is admissible: for $1\le
j\le 2\si$,
\begin{equation*}
  \|w_{\ell_j}\|_{L^\theta([t,t+\tau];L^q)}\le
  \tau^{1/\theta-1/p}\|w_{\ell_j}\|_{L^p([t,t+\tau];L^q)} \lesssim
  \tau^{1/\theta-1/p} e^{C(t+\tau)},
\end{equation*}
where we have used $(Exp)_{k-1}$. The last term of the sum is
estimated similarly. 

Finally, the term $L(w_k)$ is handled in thanks to the Gronwall lemma, and
$(Exp)_k$ follows. 
\end{proof}

The proof of Proposition~\ref{prop:growth} in the one-dimensional
cubic case follows readily. Since this case is $L^2$-subcritical, we
have $\theta<p$. Using Strichartz estimate, we infer, for $s,\tau\in \R$,
\begin{align*}
  \|u\|_{L^p([s,s+\tau];L^q)}&\le C(p)\(\|u_0\|_{L^2} +
  \|u\|^{2\si}_{L^\theta([s,s+\tau];L^q)}
  \|u\|_{L^p([s,s+\tau];L^q)}\) \\
&\le C(p)\(\|u_0\|_{L^2} +
  \tau^{2\si(1/\theta-1/p)}
  \|u\|^{2\si+1}_{L^p([s,s+\tau];L^q)}\),
\end{align*}
where we have used the conservation of mass, and $C(p)$ is independent
of $s$ and $\tau$.  Choosing $\tau$ sufficiently small, a bootstrap
argument implies that there exists $C>0$ such that
\begin{equation*}
  \|u\|_{L^p([s,s+\tau];L^q)}\le C,\quad \forall s\in \R,\ 0<\tau\le
  \tau_0. 
\end{equation*}
Again since $\theta<p$, we conclude that \eqref{eq:unifLp} is
satisfied with $f(\tau)= C \tau^{1/\theta-1/p}$. 

\subsection{Isotropic repulsive potential}
\label{sec:isorep}

We assume  $\si\ge
2/d$, and $\l\ge 0$ (defocusing nonlinearity). We show that in the
isotropic repulsive case $\Om_j=\Om\ge 0$ (a case where the 
energy $E$ defined in \eqref{eq:energy} is not a positive functional),
the evolution laws derived in \S\ref{sec:champs} show us that the
nonlinearity is negligible for large time, and there is scattering. In
this paragraph, we also assume that $\Om$ is locally Lipschitzean,
without systematically recalling this assumption. 
We 
start with the straightforward result:
\begin{lemma}\label{lem:edorep}
  Assume $\Om_j(t)\le 0$ for all $t\ge 0$. Then the solutions to
  \eqref{eq:solfondj} satisfy:
\begin{equation*}
  \nu_j(t)\ge 1,\quad \mu_j(t)\ge t,\quad \dot \nu_j(t)\ge 0,\quad
  \dot \mu_j(t)\ge 1,\quad \forall t\ge 0.
\end{equation*}
\end{lemma}
\begin{remark}
  As a consequence of this lemma, Proposition~\ref{prop:strichartz}
  remains valid with $T=\infty$, even if $\Om\le 0$ is not bounded.  
\end{remark}
We can then prove:
\begin{proposition}\label{prop:scattrepuls}
  Assume $\Om_j=\Om$ is independent of $j$, with $\Om(t)\le 0$
  for all $t\ge 0$. Let $2/d\le \si (<2/(d-2)$ if $d\ge 3$), $\l\ge 0$
  and $u_0\in \Sigma$. The solution to \eqref{eq:NLSharmo} is global
  in time, and there is scattering:
  \begin{equation*}
    \exists u_+\in \Sigma,\quad \left\lVert \Psi(t)\( u(t) -
      U(t,0)u_+\)\right\rVert_{L^2}\Tend t {+\infty} 0,
  \end{equation*}
 for any $\Psi\in \{ {\rm Id}, A_j,B_j\}$, and where $U(t,0)$
 corresponds to the free evolution \eqref{eq:linear}. 
\end{proposition}
\begin{proof}
  Since $\l\ge 0$ and $\si\ge 2/d$,
  \eqref{eq:A}, \eqref{eq:B} and Lemma~\ref{lem:edorep} yield
  \begin{equation*}
    \frac{d}{dt}\(\frac{1}{2}\|Au\|_{L^2}^2 + \frac{\l \mu^2}{\si+1}
  \|u\|_{L^{2\si+2}}^{2\si+2}\)\le 0\quad ;\quad
\frac{d}{dt}\(\frac{1}{2}\|Bu\|_{L^2}^2 + \frac{\l \nu^2}{\si+1}
  \|u\|_{L^{2\si+2}}^{2\si+2}\)\le 0. 
  \end{equation*}
We infer \emph{a priori} bounds for $Au$ and $Bu$ in $L^2$. 
Duhamel's formula reads, for $s\in\R$:
\begin{equation*}
  u(t)=U(t,s)u(s)-i\l\int_s^t U(t,\tau)\(|u|^{2\si}u(\tau)\)d\tau.
\end{equation*}
For $\Psi\in \{{\rm Id},A,B\}$, apply $\Psi$ to the above
formula:
\begin{equation*}
  \Psi(t)u(t)=U(t,s)\Psi(s)u(s)-i\l\int_s^t
  U(t,\tau)\Psi(\tau)\(|u|^{2\si}u(\tau)\)d\tau, 
\end{equation*}
where we have used the fact that $\Psi$ commutes with the linear part
of the equation. Since $\Psi$ acts on gauge invariant
nonlinearities like a derivative, we have, thanks to Strichartz
estimates:  
\begin{align*}
  \|\Psi u\|_{L^p([s,t];L^q)\cap L^\infty([s,t];L^2)}
\lesssim \|\Psi(s)u(s)\|_{L^2} +
  \|u\|^{2\si}_{L^\theta([s,t];L^q)} \|\Psi u\|_{L^p([s,t];L^q)},
\end{align*}
with the same numerology as in the proof of
Proposition~\ref{prop:Sigmaloc2}. Since $u$, $Au$ and $Bu$ belong to
$L^\infty(\R_+;L^2(\R^d))$, we have
\begin{equation}\label{eq:Lq}
  \|u(t)\|_{L^q}=\|u(t)\|_{L^{2\si+2}}\lesssim
  \frac{1}{\<t\>^{d\si/(2\si+2)}}, 
\end{equation}
where we have used the factorization formula for $A$ and $B$,
Gagliardo--Nirenberg inequality, and Lemma~\ref{lem:edorep}. We infer
$u \in L^\theta(\R_+;L^q)$:
\begin{equation*}
  \theta \frac{d\si}{2\si+2} = \frac{2d\si^2}{2-(d-2)\si}>1,
\end{equation*}
since $2d\si^2 +(d-2)\si = 2\si(d\si-1)+d\si>2$. Dividing $\R$ into a
finite number of intervals on which the $L^\theta L^q$-norm of $u$ is
small, we infer $\Psi u\in L^p(\R_+;L^q)$. Scattering follows easily:
\begin{equation*}
  U(0,t)u(t)=u_0-i\l\int_0^t U(0,s)\(|u|^{2\si}u(s)\)ds.
\end{equation*}
For $\widetilde\Psi\in \{{\rm Id},\nabla,x\}$, apply $\widetilde\Psi$
to the above formula:
\begin{align*}
 \widetilde\Psi U(0,t)u(t) &=\widetilde\Psi u_0 -i\l \int_0^t \widetilde\Psi
 U(0,s)\(|u|^{2\si}u(s)\)ds\\
&=  \widetilde\Psi u_0 -i\l \int_0^t 
 U(0,s)\Psi\(|u|^{2\si}u(s)\)ds,
\end{align*}
where $\Psi={\rm Id}$ if $\widetilde \Psi={\rm Id}$, $\Psi=-iB$ if
$\widetilde \Psi=\nabla$, and $\Psi=A$ if $\widetilde \Psi=x$,
respectively. We have 
\begin{align*}
\|\widetilde\Psi U(0,t_2)u(t_2)-   \widetilde\Psi
U(0,t_1)u(t_1)\|_{L^2} & \lesssim 
\left\lVert \int_{t_1}^{t} 
 U(0,s)\Psi\(|u|^{2\si}u(s)\)ds\right\rVert_{L^{\infty}([t_1,t_2];L^{2})}\\
& \lesssim 
\left\lVert \Psi\(|u|^{2\si}u\)\right\rVert_{L^{p'}([t_1,t_2];L^{q'})}\\
& \lesssim \|u\|^{2\si}_{L^\theta([t_1,t_2];L^q)}\| \Psi
u\|_{L^p ([t_1,t_2];L^q)}\Tend {t_1,t_2}{+\infty} 0.
\end{align*}
Therefore, $U(0,t)u(t)$ converges to some $u_+\in \Sigma$, and the
proposition follows.
\end{proof}
This result strongly suggests that the solution to the nonlinear
equation has the same behavior as the solution to the linear equation
as time goes to infinity. It should therefore not be surprising that $(Exp)_k$
is satisfied in this case. However, the delicate issue is to measure
high order Sobolev norms. To do so, we modify the argument of
Lemma~\ref{lem:growthgeneral}. We will use the operators $A$ and $B$
once, and just once in view of the discussion in the proof of
Lemma~\ref{lem:growthgeneral}. 

We have seen in the course of the
  proof of Proposition~\ref{prop:scattrepuls} that $u\in
  L^\theta(\R;L^q)$ and $\Psi u \in L^p(\R;L^q)$ for $\Psi\in \{{\rm
    Id},A,B\}$. As announced above, we modify the induction argument
  of Lemma~\ref{lem:growthgeneral}: we first apply either $A$ or $B$
  to \eqref{eq:NLSharmo}, and then apply a combination of $x^\alpha$
  and $\d_x^\beta$. We still denote by $w_\ell$ the family of
  combinations of $\alpha$ momenta and $\beta$ space derivatives, now
  applied to either $Au$ or $Bu$, with
  $|\alpha|+|\beta|=\ell-1$. Eventually, this will not alter the
  conclusion, in view of \eqref{eq:inv} and
  Lemma~\ref{lem:munuexp}. Despite this small change in the definition
  of $w_\ell$, we still have \eqref{eq:wk} for $k\ge 2$ (the case
  $k\le 1$ is of no interest, since we know that $\Psi u \in
  L^p(\R;L^q)\cap L^\infty(\R;L^2)$ for $\Psi\in \{{\rm 
    Id},A,B\}$). Resume the key estimate for $w_k$:
\begin{align*}
  \|w_k\|_{L^p([t,t+\tau];L^q)\cap L^\infty([t,t+\tau];L^2)}&\lesssim
  \|w_k(t)\|_{L^2} + \|u\|_{L^\theta([t,t+\tau];L^q)}^{2\si}
  \|w_k\|_{L^p([t,t+\tau];L^q)}\\
+ \sum_{0\le \ell_j\le
    k-1}\|w_{\ell_1}\|_{L^\theta([t,t+\tau];L^q)}&\ldots
  \|w_{\ell_{2\si}}\|_{L^\theta([t,t+\tau];L^q)}
\|w_{\ell_{2\si+1}}\|_{L^p([t,t+\tau];L^q)} 
  \\ 
& \quad + \|L(w_k)\|_{L^1([t,t+\tau];L^2)}
\end{align*}
Fixing $\tau\ll 1$ independent of $t\in \R$, the second term of the
right hand side is absorbed by the left hand side. The only difficulty
consists in analyzing the sum. We may assume that $\ell_{2\si+1}$
corresponds to the largest value of indices $\ell$.
For $1\le j\le 2\si$, if $\ell_j\le k-2$, then we simply estimate
\begin{align*}
  \|w_{\ell_j}\|_{L^\theta([t,t+\tau];L^q)}&\le
  \tau^{1/\theta}\|w_{\ell_j}\|_{L^\infty([t,t+\tau];L^q)} \lesssim
  \tau^{1/\theta}\|w_{\ell_j}\|_{L^\infty([t,t+\tau];H^1)} \\
&\lesssim
  \tau^{1/\theta} \|w_{\ell_j+1}\|_{L^\infty([t,t+\tau];L^2)}\lesssim
  e^{C(t+\tau)},
\end{align*}
where we have used $(Exp)_{k-1}$. So the only case we have to examine
is when $\ell_{2\si+1}=k-1=\ell_{j_0}$ for some $1\le j_0\le 2\si$. Note
that since $\sum \ell_j=k$, this may happen only when $k=2$. In that
case, we can assume that the term $w_{\ell_{2\si+1}}$ is of the form
$Au$ or $Bu$ (a term which is $L^p(\R;L^q)$), and estimate as above
\begin{equation*}
  \|w_{\ell_{j_0}}\|_{L^\theta([t,t+\tau];L^q)}\lesssim
  \tau^{1/\theta} \|w_{2}\|_{L^\infty([t,t+\tau];L^2)}. 
\end{equation*}
The corresponding term in the sum can therefore be absorbed by the
left hand side (like $\tt V$). In the other cases, we estimate 
$\|w_{\ell_{2\si+1}}\|_{L^p([t,t+\tau];L^q)}$ thanks to
$(Exp)_{k-1}$. Having examined all the possibilities,
we conclude that $(Exp)_k$ is satisfied. Note that
  Proposition~\ref{prop:scattrepuls} suggests that the large time
  behavior of higher (weighted) Sobolev norms of $u$ is the same as in
  the linear case, so the exponential growth is sharp in general. 

\appendix

\section{The case with no potential}
\label{sec:appen}

Consider the nonlinear Schr\"odinger equation without potential:
\begin{equation}
  \label{eq:NLS}
  i\d_t v +\frac{1}{2}\Delta v = \l |v|^{2\si}v\quad ;\quad v_{\mid
    t=0}=v_0, 
\end{equation}
with energy-subcritical or energy-critical nonlinearity, $\si \le 2/(d-2)$
if $d\ge 3$.  

\begin{lemma}\label{lem:scatt}
  Let $\si\in\N$ with $\si\ge 2$ if $d=1$, and $\si\le 2/(d-2)$ if
  $d\ge 3$. Let 
\begin{equation*}
  (p_1,q_1)= \( 2\si+2 ,\frac{2d(\si+1)}{d-2+d\si}\).
\end{equation*}
Assume that \eqref{eq:NLS} possesses a global solution $v\in
L^{p_1}(\R;W^{1,q_1}(\R^d))$. Let $k\in \N$. If $v_0\in H^k(\R^d)$, then $v\in
L^\infty(\R;H^k(\R^d))$, and more generally, $v\in
L^{p_0}(\R;W^{k,q_0}(\R^d))$ for all admissible pair $(p_0,q_0)$. 
\end{lemma}
\begin{remark}
  The assumption $\si\in \N$ is made only to simplify the
  presentation. The proof could be adapted to the case where  the
map $z\mapsto |z|^{2\si}z$ is $C^k$.  
\end{remark}
\begin{remark}
  The main assumption of the lemma states essentially that asymptotic
  completeness holds in a suitable space. We could even assume that
  the nonlinearity is $(2\si+1)$-homogeneous, and not necessarily gauge
  invariant. However, scattering is known with no size assumption on
  $v_0$ in the defocusing gauge invariant case (see below), hence our
  choice. Note 
  that in the case $d=1$, an algebraic control of the growth of
  Sobolev norms is known, regardless of gauge invariance \cite{Sta97b}. 
\end{remark}
\begin{proof}
  We remark that the pair $(p_1,q_1)$ is admissible, and 
\begin{equation*}
\frac{1}{p_1'}= \frac{2\si+1}{p_1}\quad ;\quad  \frac{1}{q_1'}=
\frac{1}{q_1}+\frac{2\si}{d\si(\si+1)}.
\end{equation*}
We prove the lemma by induction on $k$. We first prove 
\begin{equation*}
  v \in L^{p_1}(\R;W^{k,q_1})\cap L^\infty(\R;H^k). 
\end{equation*}
We start with $k=1$: applying $\nabla$ to \eqref{eq:NLS}, Strichartz
estimates on $I=[t_0,t]$ yield
\begin{align*}
  \|\nabla v\|_{L^\infty(I;L^2)\cap L^{p_1}(I;L^{q_1})} &\le C\(
  \|\nabla v(t_0)\|_{L^2} + \left\lVert 
    |v|^{2\si}\nabla v\right\rVert_{L^{p_1'}(I;L^{q_1'})}\)\\
&\le C\(
  \|v(t_0)\|_{L^2}  + \|v\|_{L^{p_1}(I;L^{d\si(\si+1)})}^{2\si}\|\nabla
  v\|_{L^{p_1}(I;L^{q_1})}\), 
\end{align*}
where we have used H\"older inequality. Notice the embedding
\begin{equation*}
  W^{1,q_1}(\R^d) \subset L^{d\si(\si+1)}(\R^d). 
\end{equation*}
In view of the assumption of the lemma, this implies $v \in
L^{p_1}(\R;L^{d\si(\si+1)})$. Therefore, we can split 
$\R$ into finitely many intervals on which
$C\|v\|_{L^{p_1}(I;L^{d\si(\si+1)})}^{2\si}\le 1/2$. On each such interval
$I$, we have
\begin{align*}
  \|\nabla v\|_{L^\infty(I;L^2)\cap L^{p_1}(I;L^{q_1})} \le 2C
  \|\nabla v(t_0)\|_{L^2} .
\end{align*}
The conclusion follows in the case $k=1$. 

Assume now that the result is known for $k\ge 1$, and that the
nonlinearity is $C^{k+1}$. Differentiating \eqref{eq:NLS} $k+1$ times
with respect to space variable, we find, for $|\alpha|=k+1$,
\begin{equation*}
  \(i\d_t +\frac{1}{2}\Delta\)\d^\alpha v = N_1(v)+N_2(v), 
\end{equation*}
with the pointwise controls
\begin{equation*}
  \lvert N_1(v)\rvert \lesssim |v|^{2\si} \lvert \d^\alpha v\rvert\quad
  ;\quad 
 \lvert N_2(v)\rvert \lesssim \sum_{|\alpha_j|\le k} \lvert \d^{\alpha_1}
v\rvert\ldots 
\lvert \d^{\alpha_{2\si+1}} v\rvert .
\end{equation*}
Strichartz estimates  on the time interval $I=[t_0,t]$ yield
\begin{align*}
  \| \d^\alpha v\|_{L^\infty(I;L^2)\cap L^{p_1}(I;L^{q_1})}&\lesssim
  \|\d^\alpha v(t_0)\|_{L^2} + 
  \sum_{j=1,2}\|N_j(v)\|_{L^{p_1'}(I;L^{q_1'})} \\
&\lesssim  \|\d^\alpha v(t_0)\|_{L^2} +
\|v\|_{L^{p_1}(I;L^{q_2})}^{2\si}\|\d^\alpha v\|_{L^{p_1}(I;L^{q_2})}\\
+\sum_{J}
\|\d^{\alpha_1}v &\|_{L^{p_1}(I;L^{q_2})}\ldots
\|\d^{\alpha_{2\si}}v\|_{L^{p_1}(I;L^{q_2})}
\|\d^{\alpha_{2\si+1}} v\|_{L^{p_1}(I;L^{q_1})} ,
\end{align*}
where we have denoted $q_2 = d\si(\si+1)$, and we have used the
ordering
\begin{equation*}
  J= \left\{ |\alpha_1|,\ldots, |\alpha_{2\si-1}|\le k-1\ ;\
 |\alpha_{2\si}|, |\alpha_{2\si+1}|\le k\quad ;\quad \sum \alpha_j= 
\alpha\right\}.
\end{equation*}
Proceeding as in the case $k=1$, we consider finitely many time
intervals on which
\begin{align*}
  \| \d^\alpha v\|_{L^\infty(I;L^2)\cap L^{p_1}(I;L^{q_1})}&\lesssim
  \|\d^\alpha v(t_0)\|_{L^2} \\ 
+\sum_{J}
\|\d^{\alpha_1}v &\|_{L^{p_1}(I;L^{q_2})}\ldots
\|\d^{\alpha_{2\si}}v\|_{L^{p_1}(I;L^{q_2})}
\|\d^{\alpha_{2\si+1}} v\|_{L^{p_1}(I;L^{q_1})} . 
\end{align*}
We use the embedding $W^{1,q_1} \subset L^{q_2}$ again, and the
induction assumption: when $\alpha_{2\si}=k$, we proceed like for the
term $N_1$ (when summing over all $\alpha$'s such that
$|\alpha|=k+1$), and in all the other cases, we deal with a
controllable source term. This yields the lemma for the pair
$(p_0,q_0)=(\infty,2)$. The estimate for general admissible pairs
follows by using Strichartz estimates again.  
\end{proof}

\begin{proposition}
  Let $\l>0$, and $\si \ge 2/d$ an integer, with $\si\le 2/(d-2)$ if $d\ge
  3$. Suppose $v_0\in \Sigma$. Let $k\in \N$, $k\ge 1$. 
  \begin{itemize}
  \item[(i)] If $v_0\in H^k(\R^d)$, then $v\in
L^{p_0}(\R;W^{k,q_0}(\R^d))$ for all admissible pair $(p_0,q_0)$. 
\item[(ii)] If in addition $x\mapsto |x|^k v_0\in L^2(\R^d)$, then
  $|x|^kv\in C(\R;L^2(\R^d))$ and for all admissible pair $(p_0,q_0)$,
  \begin{equation*}
    \forall \alpha\in \N^d, \ |\alpha|\le k,\quad \left\lVert x^\alpha
      v\right\rVert_{L^{p_0}([0,t];L^{q_0})} \lesssim \<t\>^{|\alpha|}. 
  \end{equation*}
  \end{itemize}
\end{proposition}
\begin{remark}
  In the case $\si>2/d$, $(i)$ remains valid without 
assuming $v_0\in\Sigma$. The point to notice is that one can  
  apply Lemma~\ref{lem:scatt} as in the proof below, since the
  assumptions of the lemma are known to be satisfied thanks to
  Morawetz estimates, which yield asymptotic completeness in $H^1$. In
  the case $\si=2/d$, this aspect is still an open question. 
\end{remark}
\begin{proof}
  Under our assumptions on $\l$ and $\si$, we know that there exists a unique,
  global, solution $v\in C^\infty(\R;\Sigma)$ to \eqref{eq:NLS}, with
  $v\in L^\infty(\R;H^1)$.  The
  pseudo-conformal conservation law yields
  \begin{equation*}
    \frac{d}{dt}\(\frac{1}{2}\|J(t) v\|_{L^2}^2
+ \frac{\l t^2}{\si +1}\|v\|_{L^{2\si+2}}^{2\si+2}
\)=\frac{t\l}{\si+1}(2-d\si)\|v\|_{L^{2\si+2}}^{2\si+2} ,
  \end{equation*}
where $J(t)=x+it\nabla$. The right hand side is non-positive for $t\ge
0$: this
yields an \emph{a priori} bound for $J(t) v$ in
$L^\infty(\R;L^2)$. Since 
\begin{equation}\label{eq:J}
  J(t) = it e^{i\frac{|x|^2}{2t}}\nabla\( e^{-i\frac{|x|^2}{2t}} \cdot\),
\end{equation}
Gagliardo--Nirenberg inequality yields
\begin{equation*}
  \|v\|_{L^\rho}\lesssim
  \frac{1}{|t|^{\delta}}\|v\|_{L^2}^{1-\delta}
  \|J(t)v\|_{L^2}^\delta ,\quad \text{where
  }\delta=d\(\frac{1}{2}-\frac{1}{\rho}\), \text{ and } 2\le \rho
  \le \frac{2d}{d-2}. 
\end{equation*}
We infer $v\in L^{p_1}(\R;L^{q_1})$. 
Resume the value 
\begin{equation*}
  \theta = \frac{4\si(\si+1)}{2-(d-2)\si}\quad \(\theta=\infty\text{ if
  }\si=\frac{2}{d-2}\). 
\end{equation*}
In view of the identities
\begin{equation*}
  (p,q)=\(\frac{4\si+4}{d\si},2\si+2\)\quad ;\quad
  \frac{1}{p'}=\frac{1}{p}+ \frac{2\si}{\theta}\quad ;\quad
  \frac{1}{q'}=\frac{2\si+1}{q}, 
\end{equation*}
Strichartz estimates yield
\begin{equation*}
  \|\nabla v\|_{L^p(I;L^q)\cap L^{p_1}(I;L^{q_1})}\lesssim
 1 + \|v\|_{L^\theta(I;L^{q})}^{2\si} \|\nabla v\|_{L^{p}(I;L^{q})}.
\end{equation*}
We note that $v \in L^\theta(\R;L^{q})$ (and $\|v(t)\|_{L^{q}}\to 0$
uniformly as $t\to \infty$): splitting $\R$ into
finitely many intervals, we infer $\nabla v \in L^{p_1}(\R;L^{q_1})$: the
first point of the proposition then follows from
Lemma~\ref{lem:scatt}. 
\smallbreak

The second point of the proposition is obtained by mimicking the proof
of Lemma~\ref{lem:scatt}: instead of considering $\nabla$ and its
powers, consider $J=x+it\nabla$ and its powers. In view of
\eqref{eq:J}, we can follow the same computations, since the
nonlinearity we consider is gauge invariant: $ |J|^k v\in
L^\infty(\R;L^2)$. The algebraic growth of the momenta then stems from
triangle inequality. 
\end{proof}

\subsection*{Acknowledgments} The  author wishes to thank
Clotilde Fermanian for fruitful discussions on this subject, and
Lysianne Hari for a careful reading of the manuscript.

\bibliographystyle{amsplain}
\bibliography{biblio}

\end{document}